\documentclass[a4, 11pt]{amsart}
\usepackage{amsthm,amsmath,amssymb,graphics,hyperref}
\hypersetup{pdfpagemode=FullScreen,colorlinks=true}
\usepackage{graphicx}
\usepackage{psfrag}
\usepackage{amscd}
\usepackage{enumerate}
\usepackage{varioref}
\newenvironment{thA}{\begin{trivlist}\item[]{\bf Theorem  A\ }\begin{em}}
{\end{em}\end{trivlist}}
\newenvironment{thB}{\begin{trivlist}\item[]{\bf Theorem  B\ }\begin{em}}
{\end{em}\end{trivlist}}
\newenvironment{thC}{\begin{trivlist}\item[]{\bf Theorem  C\ }\begin{em}}
{\end{em}\end{trivlist}}
\newenvironment{prop1}{\begin{trivlist}\item[]{\bf Proposition 2.2\ }\begin{em}}
{\end{em}\end{trivlist}}
\newtheorem{lemma}{Lemma}[section]
\newtheorem{thm}{Theorem}

\newtheorem{prop}[lemma]{Proposition}

\newtheorem{cl}[lemma]{Claim}

\renewcommand{\l}{l_{\rho_n}}

\newcommand{\bc}{{\mathbb C}}
\newcommand{\bh}{{\mathbb H}}
\newcommand{\bn}{{\mathbb N}}
\newcommand{\br}{{\mathbb R}}

\newcommand{\ra}{{\rightarrow}}

\newcommand{\Isom}{\mathrm{Isom}}
\newcommand{\Int}{\mathrm{Int}}
\newcommand{\sh}{\mbox{sh}}

\newcommand{\eps}{\epsilon}
\let \cal \mathcal
\title
{Convergence of  freely decomposable Kleinian groups}
\author{Inkang Kim, Cyril Lecuire and Ken'ichi Ohshika}
\date{}
\begin{document}

\maketitle

\footnotetext[3]{The first author gratefully acknowledges the
partial support of NRF Grant (NRF-2010-0024171).}

\begin{abstract}
We consider a compact orientable hyperbolic $3$-manifold with a compressible boundary. Suppose that we are given a sequence of geometrically finite hyperbolic metrics whose
conformal boundary structures at infinity diverge to a projective
lamination. We prove that if this limit projective lamination is doubly incompressible, then the sequence has  compact closure in
the deformation space. As a consequence we generalise Thurston's double limit theorem and solve his
conjecture on convergence of function groups affirmatively.
\end{abstract}
\tableofcontents
\section{Introduction}
It is one of the most important topics in the theory of Kleinian
group to study the topological structure of their deformation spaces. The quasi-conformal deformation space of a
geometrically finite Kleinian group $G$ is fairly well understood by
virtue of the work of Ahlfors, Bers, Kra, Marden and Sullivan, and
others. To put it more concretely, for a geometrically finite
Kleinian group $G$, it is known that there is a ramified covering
map from the Teichm\"{u}ller space of $\Omega_G/G$ to the
quasi-conformal deformation space of $G$, where $\Omega_G$ denotes
the region of discontinuity of $G$. On the other hand, in general it
is difficult to understand how these coordinates on Teichm\"uller
space relate to the full deformation space. In other words, it is a
challenging task to determine which divergence sequence in
Teichm\"uller space correspond to convergent sequences of Kleinian
groups.

The first example of such a  sufficient condition for convergence is
the result of  Bers in \cite{Be}, which shows that  the space of
quasi-Fuchsian groups lying on a Bers slice is relatively compact.
On the other hand, in the process of proving the uniformisation
theorem for Haken manifolds, Thurston proved the double limit
theorem for quasi-Fuchsian groups and the compactness of deformation
spaces for acylindrical manifolds, in \cite{thui} and \cite{thuii}
respectively. These are generalised to give a convergence theorem
for general freely indecomposable Kleinian groups in Ohshika
\cite{OhZ} and \cite{OhI}. The convergence in the deformation spaces
for freely decomposable groups is more complicated and is harder to
understand.

In \cite{thb}, Thurston asked how one might generalise the double
limit theorem to the setting of Schottky groups. This question was
made into a more concrete conjecture using the notion of Masur
domain, and then was generalised to function groups.
 Masur introduced in \cite{MaE} an open set in the projective  lamination space
of the boundary of a handlebody on which the mapping class group of
the handlebody acts properly discontinuously. This open set is what
we call the Masur domain nowadays. This notion is generalised by
Otal \cite{chefjeune} to the exterior boundary of a compression
body. Thurston's conjecture is paraphrased as follows: For a
sequence in Teichm{\"u}ller space converging in the Thurston
compactification to a projective lamination lying in the Masur
domain of the exterior boundary of a compression body $M$, the
corresponding sequence of convex cocompact representations in
$AH(M)$ has a convergent subsequence. Otal in \cite{conti} first
proved that Thurston's conjecture is true for rank-$2$ Schottky
space provided that the limit lamination is arational, that is, any
component of its complement is simply connected. Canary in
\cite{CanD} proved the conjecture for some special sequences in
Schottky space. Ohshika in \cite{OhA} proved the conjecture for
function groups which are isomorphic to the free products of two
surface groups under the same assumption that the limit lamination
is arational. The strongest result in this direction under the same
assumption on the limit lamination was given by Kleineidam and Souto
in \cite{boches} without any other assumption on compression bodies.
Our main result, Theorem \ref{main} yields a proof of this
conjecture of Thurston in full generality without any extra
assumption and generalises it to a slightly larger set than the
Masur domain.

We need to introduce some notions and notations to state our main theorem.
Consider a compact irreducible atoroidal $3$-manifold $M$ with boundary.
By Thurston's uniformisation theorem for atoroidal Haken manifolds, there is a
representation $\rho_0:\pi_1(M)\rightarrow \Isom(\bh^3)$ such that $\bh^3/\rho_0(\pi_1(M))$ is homeomorphic to
$\Int(M)$ by a homeomorphism that induces $\rho$.
Such a representation is said to {\it uniformise $M$}. Any quasi-conformal deformation of $\rho_0$  also uniformises $M$.
By the Ahlfors-Bers theory, when $\rho_0$ is convex-cocompact, the space $QH(\rho_0)$ of quasi-conformal deformations of $\rho_0$ up to conjugacy by  elements of $\Isom(\bh^3)$ is parametrised by the Teichm{\"u}ller space of the boundary of $M$.
More precisely, there is a (possibly ramified) covering map, called the {\it Ahlfors-Bers map} $\cal{T}(\partial M)\rightarrow
QH(\rho_0)$ whose covering transformation group is the group of isotopy classes of diffeomorphisms of $M$ which are homotopic to the identity.

The space $QH(\rho_0)$ is a subspace of the {\em deformation space} $AH(M)$.
 This deformation space $AH(M)$ is the space of discrete faithful representations $\rho:\pi_1(M)\rightarrow \Isom(\bh^3)$ up to conjugacy by elements of
$PSL_2(\bc)$.
 It is endowed with the quotient of the compact-open topology which is also called the {\em algebraic topology}.
In the main theorem, we shall consider sequences of representations given by sequences in the Teichm{\"u}ller
space
whose images under the Ahlfors-Bers map diverge in $QH(\rho_0)$ and give a sufficient condition for their convergence
in $AH(M)$.

Thurston introduced in \cite{thuiii} the notion of {\em doubly incompressible
curves}. This  can be extended to measured geodesic laminations in the following way.

 We say that a measured geodesic lamination
$\lambda\in\cal{ML}(\partial M)$ is {\em doubly incompressible} if
and only if  there exists $\eta >0$ such that  $i(\lambda,
\partial E)>\eta$ for any essential annulus or disc $E$ in $M$.
We denote by $\cal{D}(M)\subset\cal{ML}(\partial M)$ the set of doubly incompressible measured geodesic laminations and by
$\cal{PD}(M)$ its projection in the projective lamination space $\cal{PML}(\partial M)$.
It is not hard to see that $\cal{D}(M)$ contains the Masur domain (see \cite{masurext}).
Our main theorem is the following.

\begin{thm}
\label{main} Let $M$ be a compact orientable irreducible atoroidal $3$-manifold
with boundary and $\rho_0 : \pi_1(M) \ra PSL(2,\bc)$ a convex cocompact representation that uniformises $M$. Let $(m_n)$ be a
sequence in the Teichm{\"u}ller space $\mathcal{T}(\partial M)$ that
converges in the Thurston compactification to a projective measured
lamination $[\lambda]$ contained in $\cal{PD}(M)$. Let $q :
\mathcal{T}(\partial M)\rightarrow QH(\rho_0)$ be the Ahlfors-Bers
map, and suppose that $(\rho_n:\pi_1(M)\ra G_n \subset PSL(2,\bc))$
is a sequence of discrete faithful representations corresponding to
$q(m_n)$. Then  $(\rho_n)$ has an algebraically convergent
subsequence in $AH(M)$.
\end{thm}

It should be noted that our result here is closely related to the
Bers-Thurston density conjecture. This conjecture states that every
finitely generated Kleinian group is contained in the boundary of
the quasi-conformal deformation space of geometrically finite
Kleinian groups without rank-$1$ maximal parabolic subgroups. A
special case has been proved by Bromberg and Brock-Bromberg \cite{B, BB} using cone manifold deformation theory. Bromberg and Souto announced a complete prof of the general case using this approach and, in particular,  avoiding an appeal to the ending
lamination theorem (\cite{BS}).
The general case can be proved by combining the resolution of the
tameness conjecture by Agol and
 Calegari-Gabai, the ending lamination conjecture by Brock-Canary-Minsky, and
from some convergence theorems due to Thurston, Ohshika,
Kleineidam-Souto and Lecuire,  together with some topological
argument due to Ohshika and Namazi-Souto(See \cite{Ag}, \cite{CG},
\cite{thuiii}, \cite{OhA}, \cite{boches}, \cite{masurext},
\cite{Min}, \cite{OhP} and  \cite{Namazi}).


Theorem \ref{main} is a corollary of Theorem \ref{th&} below by the following argument. By Theorems of Thurston \cite{thuii}
and Canary \cite{CanD}, the convergence of $(m_n)$ to $[\lambda]$ implies that there is a
sequence of weighted multi-curves $(\lambda_n\in\cal{ML}(\partial M))$ such that
$l_{\rho_n}(\lambda_n)$ tends to $0$ and that the sequence $(\lambda_n)$ converges in
$\cal{ML}(\partial M)$ to a measured geodesic lamination whose projective class
is $[\lambda]$. Since $\lambda$ lies in $\cal{D}(M)$, Theorem \ref{main} is derived from the following
theorem, whose proof occupies the main part of this paper.

\begin{thm}
\label{th&}
Let $(\rho_n:\pi_1(M)\rightarrow \Isom(\bh^3))$ be a sequence of
convex cocompact representations that uniformise $M$ and let $(\lambda_n)\subset\cal{ML}(
\partial M)$ be a sequence of measured geodesic laminations such that
$(l_{\rho_n}(\lambda_n))$ is a bounded sequence and that $(\lambda_n)$ converges
in $\cal{ML}(
\partial M)$ to a measured geodesic lamination $\lambda\in\cal{D}(M)$. Then the sequence
$(\rho_n)$ has a compact closure in $AH(M)$, namely, any subsequence contains an
algebraically convergent subsequence.
\end{thm}

We prove Theorem \ref{th&} by contradiction using the following
arguments. We can first assume that $\lambda_n$ is a weighted simple
closed curve since the set of weighted simple closed curves is dense
in $\cal{ML}(S)$.
 Assuming that no subsequence of $(\rho_n)$ converges, we use Culler-Morgan-Shalen compactification of the character
 variety (\cite{mors1}) and ideas of Otal \cite{meister2} to construct a sequence of train tracks $\tau_n$ carrying $\lambda_n$
 and maps $f_n:\tau_n\to\bh^3/\rho_n(\pi_1(M))$ so that some parts of $\tau_n$ are mapped to long quasi-geodesic paths while the rest is mapped to relatively short paths. Then we construct a simplicial annulus joining $f_n(\lambda_n)$ to its geodesic representative $\lambda_n^*\subset M_n= \bh^3/\rho_n(\pi_1(M))$. By Gauss-Bonnet Formula, this annulus is not large enough to cover the difference in length between $f_n(\lambda_n)$ and $\lambda_n^*$. It follows that $f_n(\lambda_n)$ nearly backtracks along some long paths. We use this backtracking  to construct a sequence of annuli whose boundary converges to a sublamination of $\lambda$. This contradicts the double incompressibility of $\lambda$.\\

We state Theorem \ref{th&} under the assumption that $\rho_n$ is convex cocompact, but it should hold without this assumption, namely when $\rho_n$ is discrete and faithful. The presence of parabolics, in particular rank-$2$ cusps, would add extra technicalities to the proof, which we wished to avoid. 
Also, the presence of geometrically infinite ends would require very little change. In the same spirit, we state Theorem \ref{main} under the assumption that $\rho_n$ is convex cocompact but it should hold under the assumption that it is geometrically finite and minimally parabolic.\\

\noindent {\bf Plan of the paper:}

In Section \ref{pre}, we explain background materials and quote necessary results, some with proofs.

In Sections \ref{realisation} and \ref{track} we explain the construction of $\tau_n$ and $h_n$.
In Section \ref{realisation}, we show that for a doubly incompressible measured lamination $\lambda$ and a small minimal action of $\pi_1(M)$
on $\br$-tree $\cal T$,  each component of $\lambda$ is either realised in $\cal T$, or is carried by train tracks with arbitrarily short branches in the following sense: there is a sequence of train tracks $\theta_i$
minimally carrying the component, such that the branches of $\theta_i$ are mapped to geodesic segments in $\cal T$ with lengths going to zero.
This is proved in Lemma \ref{quivabien}.

Using this result, we construct $\tau_n$ and $f_n$ in Section
\ref{track}. By a theorem of Morgan and Shalen (\cite{mors1}), if
$(\rho_n)$ does not have a compact closure in $AH(M)$, a subsequence
of $(\rho_n)$ tends to a minimal small action of $\pi_1(M)$ on an
$\br$-tree $\cal{T}$. Let $L_{rec}$ be the union of the recurrent
leaves of the Hausdorff limit $L_\infty$ of $|\lambda_n|$ where the
$\lambda_n$ are measured laminations  with $(l_{\rho_n}(\lambda_n))$
bounded.
 We  construct a sequence of train tracks $\tau_n=\tau^1\cup\tau^2_n\cup \tau^3_n$ where $\tau^1$ carries realised components of $L_{rec}$,
$\tau_n^2$ carries non-realised components of $L_{rec}$ and
$\tau_n^3$ carries $L_\infty-L_{rec}$ with small weight with respect
to $\lambda_i$. We also construct a sequence of $\rho_n$-equivariant
maps $\hat f_n$ from the universal cover $\hat\tau_n$ of $\tau_n$ to
$\bh^3$ which map the branches of $\hat\tau_1$ to long geodesic
segments and the branches of $\hat\tau^2_n$ to comparatively short
ones. This is proved in Lemma \ref{rtt}.

In Section \ref{asum} we show that there are long parts of $f_n(\lambda_n)$ which nearly backtrack. We approximate $\lambda_n$ by weighted simple closed curves $c_n$. By hypothesis, the geodesic representatives $c_n^*\subset M_n= \bh^3/\rho_n(\pi_1(M))$ of $c_n$ have bounded lengths (taking the weight also into account).
On the other hand, by construction, $f_n(c_n)$ gets infinitely long with $n$. We construct an annulus between $f_n(c_n)$ and $c_n^*$ whose area is controlled. In this annulus, it follows from the difference in length between $f_n(c_n)$ and $c_n^*$ that, for large enough $n$, the path
$f_n(c_n)$ has a number of long segments in which it comes back
nearly parallel to itself. This provides us with some long and thin strips connecting two segments of
$f_n(c_n)$.

In Section \ref{annuli}, we explain how these strips give rise to discs or annuli whose boundaries converge in the Hausdorff topology to a geodesic lamination which
does not intersect $\lambda$ transversely.

In Section \ref{conclusion} we deduce Theorem \ref{th&} and \ref{main} from the results in the preceding sections.\\


{\bf Acknowledgements} We thank referees for numerous comments to make the paper more readable.

\section{Preliminaries}\label{pre}
\subsection{Deformation space}
Let $G$ be a finitely generated torsion-free Kleinian group, namely
a (torsion-free and finitely generated) discrete subgroup of
$\Isom^+(\bh^3)$. The group $G$ is {\em convex cocompact} if there is a compact subset $C\subset\bh^3/G$ that contains all the closed geodesics. Denote by $\Omega_G$ the domain of discontinuity for the action of $G$ on $\hat\bc$. The group $G$ is convex cocompact if and only $\bh^3/G\sqcup\Omega_G/G$ is compact.

Let $M$ be a compact orientable irreducible atoroidal manifold. Let $AH(M)$
denote the set of faithful discrete representations from $\pi_1(M)$
to $PSL(2, \bc)$ modulo conjugacy. We endow $AH(M)$ with the
topology induced from the representation space. The subspace of
$AH(M)$ consisting of convex cocompact representations is denoted by $CC(M)$. This space $CC(M)$ is not empty if and only if $\partial M$ contains no tori and it may contain several connected components.
The component consisting of representations $\rho$ for which there
is a homeomorphism from $\mathrm{Int}(M)$ to $\bh^3/\rho(\pi_1(M))$ that induces $\rho$ is denoted by $CC_0(M)$. For such representations, the homeomorphism $\mathrm{Int}(M)\to\bh^3/\rho(\pi_1(M))$ extends to a homeomorphism $M\to \bh^3/\rho(\pi_1(M))\sqcup\Omega_\rho/\rho(\pi_1(M))$. This produces a natural identification of $\Omega_\rho/\rho(\pi_1(M))$ with $\partial M$. Notice that the homeomorphism $M\to \bh^3/\rho(\pi_1(M))\sqcup\Omega_\rho/\rho(\pi_1(M))$ is well defined up to composition by an element of $\mathrm{Mod}_0(M)$ (the set of diffeomorphisms of $M$ that are homotopic to the identity map). 

For a Kleinian group $G$, if there is a quasi-conformal automorphism
$f$ of $S^2_\infty$ such that $f Gf ^{-1}$ is again a Kleinian
group, then this group $f G f^{-1}$ is said to be a {\em
quasi-conformal deformation} of $G$. We denote by $QH(G)$ the space of quasi conformal deformations of $G$ up to conjugacy. By the theory
of Ahlfors-Bers, there is a ramified covering map from
$\mathcal{T}(\Omega_G/G)$ to the space of quasi-conformal
deformations of $G$ modulo conjugacy. 

Given a representation $\rho_0\in CC_0(M)$, we have $QH( \rho_0(\pi_1(M))=CC_0(M)$. As we have seen, for any $\rho\in CC_0(M)$, there is a natural identification of $\Omega_\rho/\rho(\pi_1(M))$ with $\partial M$. Thus the  theory
of Ahlfors-Bers provides us with a covering map $\mathcal{T}(
\partial M) \rightarrow CC_0(M)$ which we call the {\em Ahlfors-Bers map}.

\subsection{$\br$-trees}\label{tree}
An $\br$-tree $\cal{T}$ is a geodesic metric space in which any two points $x,y$ can be
joined by a unique simple arc. Let $G$ be a group acting by isometries on an $\br$-tree
$\cal{T}$. The action is {\em minimal} if there is no proper invariant subtree and {\em
small} if the stabilizer of any non-degenerate arc is virtually Abelian.

Morgan and Shalen \cite{mors1} made use of $\br$-trees to compactify
deformation spaces. They used algebraic methods involving
valuations, while the same result has been obtained by Paulin
\cite{paulin} and Bestvina \cite{best} using a more geometrical
approach. In this paper we shall adopt the point of view of
Kapovich-Leeb \cite{kl} (see also \cite[chapters 9 and 10]{kapo}).
Let $(\rho_n\subset GF_0(M,P))$ be a sequence of representations
such that no subsequence of $(\rho_n)$ converges algebraically. Let
$\Gamma \subset\pi_1(M)$ be a  set of generators and let $\tilde
x_n\subset\bh^3$ be a point realising the minimum $\eps^{-1}_n$ on
$\bh^3$ of the function $\max\{d(\tilde x,\rho_n(a)(\tilde x)),a\in
\Gamma\}$ (see for example \cite[Lemma 6.5]{paulin} for the existence of such a
point). Since no subsequence of $(\rho_n)$ converges algebraically,
$(\eps_n^{-1})$ tends to $\infty$. Choose a non-principal
ultra-filter $\omega$ and denote by $\eps_n\bh^3$ the hyperbolic
space $\bh^3$ with the hyperbolic metric rescaled by $\eps_n$. The
ultra-limit $(X_\omega,x)=\omega-\lim (\eps_n\bh^3,x_n)$ of the
sequence of rescaled spaces is defined as follows. Let
$\Pi_{n}(\eps_n\bh^3)$ be the infinite product of the spaces
$(\eps_n\bh^3)$. We define a function $d_\omega$ on
$\Pi_{n}(\eps_n\bh^3)$ by setting $$ d_\omega(y,z)=\omega-\lim
d_{\eps_n\bh^3}(\tilde y_n,\tilde z_n)$$ for any two points
$y=(\tilde y_n)$ and $z=(\tilde z_n)$ lying in
$\Pi_{n}(\eps_n\bh^3)$.

This function $d_\omega$ is a pseudo-distance in
$\Pi_{n}(\eps_n\bh^3)$ with values in $[0,\infty]$ and we set
$(X_\omega,d_\omega)=(\Pi_{n}(\eps_n\bh^3),d_\omega)/\sim$ where we
identify points with zero $d_\omega$-distance. Let $x=(\tilde x_n)$
denote the sequence of points $\tilde x_n$ defined above. The metric
space $(X_\omega,x)$ is the set of points of $(X_\omega)$ with a
finite distance from $x$. This metric space is an $\br$-tree (cf.
\cite{kl}). The action of $\rho_n(\pi_1(M))$ on $\eps_n\bh^3$ gives
rise to an action of $\pi_1(M)$ on $(X_\omega,x)$ by isometries.
This action is small (cf. \cite{kl}). Let $\cal{T}$ be the minimal
invariant subtree of $X_\omega$ under this action. We say that
$(\rho_n)$ tends to the action of $\pi_1(M)$ on $\cal{T}$ with
respect to $\omega$. For $c\in\pi_1(M)$ let us denote by
$\delta_\cal{T}(c)$ the minimal translation distance of $c$ on
$\cal{T}$. Then we have $\delta_\cal{T}(c)=\omega-\lim \eps_n
l_{\rho_n}(c)$, where $l_{\rho_n}(c)$ is the length in
$\bh^3/\rho_n(\pi_1(M))$ of the closed geodesic in the free homotopy
class of $c$.

\subsection{Geodesic laminations}
\indent A {\em geodesic lamination} $L$ on a complete hyperbolic
surface $S$ is a compact set which is a disjoint union of complete
embedded geodesics called {\em leaves}. It is known that this
definition is independent of a chosen hyperbolic metric on $S$.
The details can be found in  \cite{meister2}.
For a connected geodesic lamination
$L$ which is not a simple closed curve we denote by $\bar S(L)$ the
smallest subsurface of $S$ with compact geodesic boundary containing
$L$. Inside $\bar S(L)$ there are finitely many closed geodesics
(including the components of $
\partial\bar S(L)$) disjoint from $L$.
These closed geodesics do not intersect each other
(cf. \cite{espoir} page 99) and we denote by $
\partial'\bar S(L)\supset
\partial\bar S(L)$ their disjoint union.
For example if a component of $\bar S(L)\setminus L$ is an annulus, its core curve
is contained in $\partial'\bar S(L)$ but not in $\partial\bar S(L)$.
Removing a small open tubular neighbourhood of $\partial'\bar S(L)$ from
$\bar S(L)$ we get a compact surface $S(L)$. We call $S(L)$ the {\it
surface embraced} by the geodesic lamination $L$ and $\partial'\bar
S(L)$ the {\it effective boundary} of $S(L)$. If $L$ is a simple
closed curve, we define $S(L)$ to be an annular neighbourhood of $L$
and we take $\partial'\bar S(L)=L$. When $L$ is not connected,
$S(L)$ is the disjoint union of the surfaces embraced by the
connected components of $L$.

We say that a geodesic measured lamination  $L$ {\em crosses} another geodesic lamination  $L'$  if {\em at least} one leaf of $L$ intersects
a leaf of $L'$ transversely.

 A {\em measured geodesic lamination} $\lambda$ is a geodesic lamination
$|\lambda|$ together with a transverse measure of full support. We
denote by $\cal{ML}(S)$ the space of measured geodesic laminations
on $S$ endowed with the weak-$*$ topology. To simplify the
notations, we write $\cal{ML}(\partial M)$ instead of
$\cal{ML}(\partial_{\chi<0} M)$ for a compact $3$-manifold $M$ with
boundary. The projective lamination space $\cal{PML}(\partial M)$ is
defined to be $(\cal{ML}(\partial M) -\{0\})/\br_+^*$ where $0$
stands for the measured lamination with empty support. It should be
noted that $\cal{ML}(\partial M)$ contains measured laminations
whose restriction to some component of $\partial M$ is empty. The
Teichm{\"u}ller space $\cal{T}(\partial M)$ denotes similarly
$\cal{T}(\partial_{\chi<0}M)$. The boundary of the Thurston
compactification of $\cal{T}(\partial M)$ is equal to
$\cal{PML}(\partial M)$.
\subsection{Some notations}\label{notations}
When $\lambda$ is a  measured geodesic lamination, we denote by
$|\lambda|$ the support of $\lambda$. For an arc $k$ whose
intersections with $|\lambda|$ are transverse, we will denote by
$\int_k d\lambda$ the $\lambda$-measure of $k$.

Let $(u_n)$ and $(v_n)$ be two sequences of non-negative real
numbers. We say that $u_n$ is $o(v_n)$ if for any $\eps$ there is
$N(\eps)$ such that for $n\geq N(\eps)$, we have $u_n\leq\eps v_n$.
We also write $u_n=o(v_n)$.

We will say that $u_n$ is $O(v_n)$ if there are $K,N>0$ such that
for $n\geq N$, we have $u_n\leq K v_n$.

We say that $u_n$ is $\Theta(v_n)$ if $u_n$ is $O(v_n)$ and $v_n$ is
$O(u_n)$.

\subsection{Train tracks and their realisations}

A {\em train track} $\tau$ in a hyperbolic surface $S$ is a  union
of finitely many rectangles with a distinguished pair of vertical
opposite sides. These rectangles meet each other only along non-degenerate
segments contained in their vertical sides in such a way that every point on a
vertical side of a rectangle lies in at least one other rectangle.
The rectangles are called {\em branches}, and they are foliated by
vertical segments called {\em ties}. A connected component of the
intersection of the branches is called a {\em switch}. The branches
are also foliated by horizontal segments, and a smooth arc which is
a union of horizontal segments is called a {\em rail} or a {\em
train route}. A geodesic lamination is {\em carried} by $\tau$ if it
is isotopic to one which lies in $\tau$ in such a way that the leaves
are transverse to the ties (see \cite{Bouts} or \cite{meister2} for
more details about train tracks).

When $\tau$ is a train track and $\lambda$ is a measured geodesic
lamination whose support is carried by $\tau$, we say that $\lambda$
is carried by $\tau$. For a branch $b$ of $\tau$, we define the
number $\lambda(b)$ to be the $\lambda$-measure of a tie of $b$.
This number does not depend on the choice of a tie in $b$.

Consider an action of $\pi_1(S)$ on an $\br$-tree $\cal{T}$ by
isometries. A measured lamination $\lambda$ is said to be {\em
realised} in $\cal{T}$ if there is a $\pi_1(S)$-equivariant map
$\phi: \bh^2 \ra {\cal T}$ such that the restriction of $\phi$ to
any lift of a leaf of $\lambda$ in $\bh^2$ is a complete geodesic in
the tree, and identifying this geodesic with $\mathbb{R}$ gives
a weakly monotonous unbounded function on the leaf.
A train track
$\tau\in S$ is said to be realised in $\cal{T}$ if there is an
equivariant map $\phi: \bh^2 \ra {\cal T}$ which maps each lift of a
branch of $\tau$ to a non-trivial geodesic segment on ${\cal T}$ in
such a way that each rail is mapped injectively, and that each tie
collapses to a point. By \cite{meister2}, $\lambda$ is realised in
$\cal{T}$ if and only if $\lambda$ is carried by a train track
$\tau$ which is realised in $\cal{T}$.

We say that a measured lamination $\lambda$ is {\em collapsed} by
$\phi: \bh^2 \ra {\cal T}$ when there is a train track $\tau$
carrying $\lambda$ such that every component of $\tilde{\tau}$, the
preimage of $\tau$ in $\bh^2$, is mapped to a point by $\phi$. It is
straightforward that there exists an equivariant map collapsing
$\lambda$ if and only if the action of $i_*(\pi_1(S(\lambda))$ on
$\cal{T}$ has a global fixed point.

For a  measured lamination $\mu$ on a surface $S$, if it does not have atoms, the semidistance on $\bh^2$ induced by integrating the transverse
measure $\tilde \mu$ along paths is continuous with respect to the usual topology of $\bh^2$. By replacing closed leaves by annuli foliated
by parallel closed curves, we obtain a measured partial foliation $\cal F_\mu$. The quotient of $\bh^2$ under the semi-distance induced by
$\cal F_\mu$ gives rise to a dual tree $\cal T_\mu$ with the projection $\pi:\bh^2\ra \cal T_\mu$.

A morphism $\phi:\cal T\ra \cal T'$ between $\br$-trees is a map with the property that every non-degenerate arc $[p,q]\subset \cal T$
contains a non-degenerate subarc $[p,r]\subset [p,q]$ which is mapped isometrically onto $\phi[p,r]\subset \cal T'$. A morphism from
$\cal T_\mu$ to $\cal T$ is said to fold only at complementary regions if the only folding points are projections of complementary
regions of $\tilde{ \cal F_\mu}$, where $p$ is a folding point if $[p, q],[p, q']\subset \cal T_\mu, [p,q]\cap [p,q']=\{p\}$ is mapped to the same segment in $\cal T$.
The following theorems \cite{chefm} will be useful to us later.
\begin{thA}{\rm [Morgan-Otal]}\label{morphism} Let $(\alpha_1,\cdots,\alpha_k)$ be a collection of simple closed curves which
define a pants decomposition of a surface $S$ and let $\pi_1(S)$ act on an $\br$-tree $\cal T$.
Then there is a measured lamination $\mu$ on $S$ and an equivariant morphism
$$\phi:\cal T_\mu\rightarrow \cal T$$ with $\delta_{\cal T}(\alpha_i)=\delta_{\cal T_\mu}(\alpha_i)=i(\mu,\alpha_i)$ for all $i$.
Moreover $\phi$ folds only at complementary regions.
\end{thA}

The above theorem is generalised by Skora \cite{skora}.
\begin{thB}{\rm [Skora]} Suppose the action of $\pi_1(S)$  on an $\br$-tree $\cal T$ is minimal and small such that the action of each element
representing $\partial S$ has a fixed point in $\cal T$. Then there is a unique measured
lamination $\mu$ and an equivariant isometry
$$\phi:\cal T_\mu\rightarrow \cal T.$$
\end{thB}

\subsection{Compactification of $\tilde M$}\label{compactification}

We denote by $\tilde M$ the universal covering of $M$ and by
$p:\tilde M\rightarrow M$ the covering projection. We compactify
$\tilde M$ in the following way: endow $M$ with a geometrically
finite hyperbolic metric $\sigma$ with minimal parabolics which does not correspond to a Fuchsian representation, and let us
denote by $N(\sigma)^{\rm thick}$ the complement in the convex core
$N(\sigma)$ of $\epsilon$-thin neighbourhood of the cusps of
$\sigma$ for some $\epsilon$ smaller than the Margulis constant. Let
us choose an isometry between the interior of $\tilde M$ and $\bh^3$. Now we can consider $\tilde
N(\sigma)^{\rm thick}$ as a closed subset of $\bh^3$. Since $\sigma$
is geometrically finite, there is a natural homeomorphism between
$M$ and $N(\sigma)^{\rm thick}$. Therefore we can regard $\tilde M$
as a closed subset of $\bh^3$. The compactification
$\overline{\tilde M}$ of $\tilde M$ is the closure of this closed
subset in the usual compactification of $\bh^3$ by the unit ball. If
we replace $\sigma$ by another geometrically finite metric $\sigma'$
with minimal parabolics, it follows from results of \cite{floyd}
that we get a compactification which is homeomorphic to the one
obtained with $\sigma$. Therefore this definition is independent of
the metric we chose. We call this the Floyd-Gromov compactification
of $\tilde M$.
To denote the ideal boundary in the Floyd-Gromov compactification, we use the symbol $\partial_\infty \tilde M$.

\indent A {\em meridian} is a simple closed curve $c\subset
\partial M$ which bounds a disc in $M$ but not on $\partial M$. A compact surface
$\Sigma\subset
\partial M$ is incompressible if it contains no meridians. When we consider the
closure of a lift of an incompressible surface in $\overline{\tilde
M}$, we have the following:

\begin{lemma} \label{disque} Let $\Sigma\subset
\partial_{\chi<0} M$ be a compact connected incompressible surface which does not contain any
essential closed curve homotopic into $\partial_{\chi=0} M$.
Let $\hat\Sigma$ be the
universal covering of $\Sigma$, which is completed
in the usual way to a closed disc $\overline{\hat\Sigma}$. Then any lift of $\Sigma$ to $
\partial\tilde M$ is a disc whose closure in $\overline{\tilde M}$ is homeomorphic to
$\overline{\hat\Sigma}$ in an equivariant way.
\end{lemma}

\begin{proof}
This lemma was proved in \cite{chefjeune}. See also Lemme 2.4 of \cite{espoir}.
\end{proof}

Let $\Sigma\subset
\partial M$ be a compact (possibly disconnected) incompressible surface. Johannson and Jaco-Shalen defined a
characteristic submanifold $W$ relative to $\Sigma$ (cf. \cite{jojo}
and \cite{jacs}).
We say that an $I$-bundle  $N\subset M$ is essential in
$(M,\Sigma)$ if
\begin{enumerate}
\item  $\pi_1(N)$ injects to $\pi_1(M)$ by the homomorphism induced from the inclusion,
\item $N\cap \partial M\subset \Sigma$ is the associated $\partial I$-bundle, and
\item $N$ cannot be homotoped into $\Sigma$ by a homotopy fixing $N\cap\Sigma$.
\end{enumerate}
Similarly a solid torus $D^2\times S^1$ is essential in $(M,\Sigma)$
if it satisfies
 the first and the third conditions above.
%

When we are considering an atoroidal manifold $M$ whose boundary contains no torus,  a characteristic
submanifold is a disjoint union of essential (possibly twisted) $I$-bundles over
compact incompressible surfaces and essential solid tori.
A disjoint
union $W$ of essential (possibly twisted) $I$-bundles and essential solid tori is said to be a
characteristic submanifold if and only if it has the following two
properties:\smallskip

\begin{itemize}
\item any essential (possibly twisted) $I$-bundle and any essential solid torus in $(M,\Sigma)$ can be homotoped in $W$;
\item no connected component of $W$ can be homotoped into another connected component of
$W$.\smallskip
\end{itemize}

By \cite{jojo} and \cite{jacs}, if $W$ and $W'$ are two  characteristic submanifolds
relative to $\Sigma$, then there is a diffeomorphism $\psi:M\rightarrow M$ isotopic to the
identity relative to $\partial M-\Sigma$ such that $\psi(W)=W'$ and that $\psi
(W\cap\Sigma)=W'\cap\Sigma$.

Such a characteristic submanifold can be found by looking only at $\overline{\tilde M}-\tilde
M$.

\begin{prop}[\cite{espoir}, \textsection $2$, paragraphs after Lemme 2.7] \label{carac} Let $\Sigma$ and $\Sigma'\subset
\partial_{\chi<0} M$ be two compact, connected, incompressible surfaces which are disjoint
or equal and do not contain any essential closed curve which can be homotoped into $
\partial_{\chi=0} M$. Let $\tilde\Sigma\subset
\partial\tilde M$ (resp. $\tilde\Sigma'$) be a connected component of the preimage of
$\Sigma$ (resp. $\Sigma'$) and let $\Gamma\subset\rho(\pi_1(M))$ (resp. $\Gamma'$) be the
stabiliser of $\tilde\Sigma$ (resp. $\Gamma'$).

Then $\overline{\tilde\Sigma} \cap\overline{\tilde\Sigma}'$ is either empty or
equal to the limit set of $\Gamma \cap\Gamma'$.

In the latter case, if $\Gamma \cap\Gamma'$ is not cyclic, then it is the fundamental
group of a (possibly twisted) $I$-bundle which is a connected component of a characteristic submanifold of
$(M,\Sigma \cup\Sigma')$. If $\Gamma \cap\Gamma'$ is cyclic, then it is a finite index
subgroup
of a solid torus which is a connected component of a characteristic submanifold of
$(M,\Sigma \cup\Sigma')$.
\end{prop}
See Appendix for a brief proof of Proposition \ref{carac}.

\subsection{Geodesic laminations on compressible surfaces}
Let $M$ be a compact $3$-manifold with boundary, and let $c\subset\partial M$ be a simple
closed curve. A {\em $c$-wave} is a simple arc $k$, with $k\cap c=
\partial k$ such that there is an arc $\kappa$ in $c$ with the simple closed curve
$k\cup\kappa$ bounding a compressing disc in $M$.
In some literature, a $c$-wave is allowed to intersect $c$ in its interior.
A simple innermost argument shows that if there is a $c$-wave in this
generalised sense, there is one in our sense.

Let $L$ be a geodesic lamination on $
\partial_{\chi<0} M$, and let $c\subset
\partial M$ be a multi-curve.
In the following, we always assume that simple closed curves or
multi-curves are geodesics for a fixed reference hyperbolic metric, hence there are no inessential
intersections between them or with geodesic laminations. We say that
$L$ is in {\em tight position} with respect to $c$ if $L$ contains
no $c$-waves and if {\em every} leaf of $L$ intersects $c$
transversely.

Following \cite[Theorem 1.6]{chefjeune}, we use cut-and-paste operations to construct a meridian $m$ such that a given geodesic lamination contains no $m$-waves.

\begin{cl}   \label{limit}
Let $F\subset
\partial M$ be a compressible compact surface, and let $\beta\subset F$ be a measured geodesic
lamination. Then either $\beta$ intersects transversely a meridian
$m$ and contains no $m$-waves, or there is a sequence of meridians
$(m_i\subset F)$ converging in the Hausdorff topology to a geodesic
lamination which does not cross
$\beta$.
\end{cl}
\begin{proof}
If $\beta$ intersects no meridians transversely, then, since $F$ is compressible,
there is a meridian $m\in F$ such that $i(\beta,m)=0$. Setting $m_i=m$ for every $i$,
we get
the conclusion.

Now we assume that $\beta$ intersects a meridian $m$ transversely.
If $\beta$ contains an $m$-wave $k$, let us ``cut $m$ along $k$'' to get a new meridian $m_1$
: let $\kappa$ be the closure of a connected component of $m-\partial k$ such that
$\int_{\kappa}d \beta\leq \frac{1}{2} i(\beta,m)$,  and $m_1$ the simple closed geodesic
in the
free homotopy class of $k\cup\kappa$.
We have $i(\beta , m_1)\leq\frac{1}{2} i(\beta,m)$
and $m_1$ is a meridian. If $\beta$ contains no $m_1$-waves, we are done.
If $\beta$ contains an $m_1$-wave $k_1$, then we cut $m_1$ along $k_1$ as above to get a new
meridian $m_2$ with $i(m_2,\beta)\leq\frac{1}{4} i(m,\beta)$.
Repeating this operation,
we get either a meridian $m'$ such that $\beta$  contains no $m'$-waves or a
sequence of meridians $(m_i)$ such that $i(m_i,\beta)\longrightarrow 0$.
In the latter case, we can extract a subsequence such that $m_i$ converges in the Hausdorff topology to some
geodesic lamination $H$.
Since $i(m_i,\beta)\longrightarrow 0$, we see that $H$ does not cross $\beta$.
\end{proof}

Much of the rest of this section parallels the argument in
\cite{espoir}, and we refer the readers there for more details.

We shall now introduce homoclinic
leaves. They have already proved to be a useful tool in the study of laminations on the compressible boundary of a hyperbolic $3$-manifold (see for example \cite{chefjeune}, \cite{boches} and \cite{espoir}). They are also related to Culler-Morgan-Shalen compactification of Character Varieties of freely decomposable Kleinian groups (see \cite{boches} and \cite{espoir}). They shall play an important role in different places in this paper.

Let $l$ be a simple geodesic on $
\partial_{\chi<0} M$. Such a simple geodesic $l$  is said to be {\em homoclinic} if a lift $\tilde l$
of $l$ to
the universal cover $\tilde M$ of $M$ contains two sequences of point $(\tilde x_n)$ and
$(\tilde y_n)$ such that the distance between $\tilde x_n$ and $\tilde y_n$ in $\tilde M$
is uniformly bounded whereas their distance measured on $\tilde l$ tends to $\infty$. By
Lemma \ref{disque}, an incompressible surface cannot contain a homoclinic geodesic.

Homoclinic leaves appear naturally in Hausdorff limits of sequences of meridians.
This is illustrated in the following criterion of Casson whose proof can be found in
\cite{chefjeune} and \cite[Theorem B.1]{espoir}.

\begin{lemma}
\label{casson}
Let $(m_n\subset
\partial M)$ be a sequence of meridians which converges to a geodesic lamination $H$ in the
Hausdorff topology. Then $H$ contains a homoclinic leaf.
\end{lemma}

\indent
A {\em simple half-geodesic} is an embedded half-line in $\partial M$ whose image is locally geodesic for some hyperbolic
metric on
$\partial_{\chi<0} M$.
 Let $\tilde l_+\subset\partial\tilde M$ be a half-geodesic and let $\Bar{\tilde l}_+$ be its closure in the Floyd-Gromov compactification of $\tilde M$.
We  say that $\tilde l_+$ {\em has a well-defined endpoint} if $\Bar{\tilde l}_+-\tilde l_+$ contains one point.
We  say that a geodesic $\tilde l\subset\partial\tilde M$ has two well-defined endpoints if $\tilde l$ contains two disjoint half geodesics each having a well-defined
 endpoint.
 Notice that we allow the two endpoints to be the same.\\

As in the introduction, we say that a measured geodesic lamination
$\lambda\in\cal{ML}(\partial M)$ is {\em doubly incompressible} if and only if
there exists $\eta>0$ such that  $i(\lambda,
\partial E)>\eta$ for every essential annulus, M{\"o}bius band or disc $E$ in $M$.\smallskip\\
\indent
We denote by $\cal{D}(M)\subset\cal{ML}(
\partial M)$ the set of doubly incompressible measured geodesic laminations and by
$\cal{PD}(M)$ its projection in the space $\cal{PML}(
\partial M)$ of projective measured laminations.\\

\indent Some properties of this set ${\cal D}(M)$ are discussed in
\cite{masurext}. One can deduce the following from \cite{espoir}.

\begin{lemma}
\label{saintptot}
Let $\lambda\in{\cal D}(M)$ be a measured geodesic
lamination and $l_+, l_-\subset \partial M$ two simple half-geodesics which do not intersect $|\lambda|$
transversely.
 Then any lift of $l_+$ (resp. $l_-$) to $\tilde M$ has a well-defined endpoint in $\partial_\infty \tilde M$.

Furthermore, if a lift of $l_+$ has the same endpoint as a lift of $l_-$, then they are asymptotic on $\partial \tilde M$.
\end{lemma}

\begin{proof}
The relation between doubly incompressible laminations and bending laminations is explained in \cite[Lemma 3.5]{masurext}. Knowing this relation, the proof of Lemma \ref{saintptot} can be deduced from \cite{espoir} as follows.

The first property, namely that any lift of $l_+$ (resp. $l_-$) to $\tilde M$ has a well-defined endpoint in $\partial_\infty \tilde M$, can be deduced from the proofs of
 \cite[Lemme 3.1]{espoir} and \cite[Lemme 3.3]{espoir}.

The proof of the second property, namely that if a lift of $l_+$ has the same endpoint as a lift of $l_-$ then they
are asymptotic on $\partial \tilde M$,
can be found in Lemme C5 and more specifically in the paragraph after  \cite[Affirmation C3]{espoir}.
\end{proof}

From \cite{masurext}, we also get the following:

\begin{lemma} \label{anho}
Let $\lambda \in {\cal D}(M)$ be a measured geodesic lamination and $h$ a homoclinic simple geodesic. Then the
support $|\lambda|$ of $\lambda$ crosses $h$.
\end{lemma}

\begin{proof} When $M$ is not a genus-$2$ handlebody, this is \cite[Lemma 3.6]{masurext}.
The case when $M$ is a genus-$2$ handlebody is discussed in \cite{masurext} in the remark
following \cite[Lemma 3.6]{masurext}.
\end{proof}

By the same argument as the proof of Claim \ref{limit}, we get the following.
\begin{lemma}\label{tightposition} Let $\lambda$ be a measured geodesic lamination in $\cal D(M)$, and let $S\subset
\partial
M$ be a compressible surface. Then there is a meridian $m$ in $S$ such that $S$ contains no $m$-waves disjoint from $|\lambda|$.
\end{lemma}
\begin{proof}
Take any meridian $m$ in $S$.
If there is an $m$-wave in $S$ which is disjoint from $|\lambda|$, then by the same argument as the proof of Claim \ref{limit}, we get a sequence of meridians $(m_i)$ on $S$ with $i(m_i, \lambda) \longrightarrow 0$.
This contradicts the assumption that $\lambda$ is doubly incompressible.
\end{proof}



\section{Realisations of doubly incompressible laminations}\label{realisation}

Let $\pi_1(M)\curvearrowright\cal{T}$ be a small minimal action of
$\pi_1(M)$ on an $\br$-tree for a compact irreducible atoroidal
$3$-manifold $M$. Let $S$ be a connected component of
$\partial_{\chi<0} M$. Using the map
$i_*:\pi_1(S)\rightarrow\pi_1(M)$ induced by the inclusion, we get
an action of $\pi_1(S)$ on ${\cal T}$. Therefore, if
$\lambda\in\cal{ML}(S)$ is a measured geodesic lamination, it makes
sense to ask whether or not it is realised in $\cal{T}$. In this
section, we shall discuss this question for the connected components
of a measured lamination lying in $\cal{D}(M)$.

If a lamination is
collapsed to a point by an equivariant map to a tree, then it is not
realisable.
Although the converse is not true,  one may expect that non-realisability of a
minimal component implies that the component is carried by a train
track of  bounded complexity whose edges are mapped to arbitrarily short arcs in
the $\br$-tree. In the following lemma, we shall show that this expectation is fulfilled for a doubly incompressible lamination. The proof will occupy the rest of the section.

\begin{lemma}  \label{quivabien} Let $\pi_1(M)\curvearrowright\cal{T}$ be a
small minimal action of $\pi_1(M)$ on an $\br$-tree. Let $\lambda\in\cal{D}(M)$
be a measured geodesic lamination, and $\alpha$ a minimal sublamination of $\lambda$.
Then one of the following holds:

\begin{itemize}
\item the measured lamination $\alpha$ is realised in $\cal{T}$;
\item there is a sequence of train tracks $\theta_i$ each of which minimally
carries $\alpha$ and has the following properties:
\begin{itemize}
\item $\theta_i$ has only one switch $\kappa_i$
and $\kappa_i\supset\kappa_{i+1}$ for every $i$.
\item There are a sequence
$\eta_i\longrightarrow 0$, and a sequence of $\pi_1(S)$-equivariant maps
$\phi_i:\bh^2\rightarrow\cal{T}$ such that $\phi_i$ maps every branch of the
preimage of $\theta_i$ to a geodesic segment (which may be a point) with length
smaller than $\eta_i$.
\item $\phi_i(\kappa_i)=x$ does not depend on $i$.
\end{itemize}
\end{itemize} Furthermore there is at least one component of $\lambda$ which
is realised in $\cal{T}$.
\end{lemma}
Note that, since $\theta_i$ has only one switch, the number of its branches is uniformly bounded.
\begin{proof}
When $\alpha$ is collapsed in $\cal{T}$, it is easy to see that we are in the
second situation.
The difficult case here is when $\alpha$ is neither collapsed nor realised in $\cal{T}$.

When $\alpha$ is a simple closed curve, either $\alpha$ is collapsed (when $\delta_\cal{T}(\alpha)=0$) or $\alpha$ is realised (when $\delta_\cal{T}(\alpha) > 0$).

From now on, we assume that $\alpha$ is not a simple closed curve. As we will see later, the proof is easier when $\partial M$ is incompressible. Bearing that in mind, we first get rid of the meridians that are disjoint from $S(\lambda)$. We cut $M$ along
a maximal family of compressing discs disjoint from $S(\alpha)$. We denote by $N$ the
connected component of the resulting manifold that contains $\alpha$ on its boundary.
By construction, the surface $\partial N-S(\alpha)$ is incompressible in $N$.

Using the homomorphism $i_*:\pi_1(N)\rightarrow \pi_1(M)$ induced by
the inclusion, we view $\pi_1(N)$ as a subgroup of $\pi_1(M)$. Thus
we get a small action of $\pi_1(N)$ on $\cal{T}$. Let $\cal{T}_N$ be
the minimal subtree of $\cal{T}$ that is invariant under the action
of $\pi_1(N)$ regarded as a subgroup of $\pi_1(M)$.

We denote by $F$ the component of $\partial N$ that contains
$\alpha$. Since $\lambda$ lies in $\cal{D}(M)$, no component of
$\partial S(\alpha)$ bounds a disc in $M$. It follows that the
measured lamination $\alpha$ can be regarded as an element of
$\cal{ML}(F)$. Using the map $i_*:F\rightarrow N$ induced by the
inclusion, we get an action of $\pi_1(F)$ on $\cal{T}_N$ (which is
not small when $\partial N$ is compressible). Let $\cal{T}_F$ be the
minimal subtree of $\cal{T}_N$ that is invariant under the action of
$\pi_1(F)$. If $\cal{T}_F$ is trivial, then $\alpha$ is collapsed.
From now on we assume that $\cal{T}_F$ is not trivial.

As is often the case for $3$-manifolds, the case when $\partial N$ is incompressible is easier and will give us important insights for the general case. Since we are only interested in $\alpha$, we just need to assume that $F$ is incompressible.\\

 \underline{First case} : $F$ is incompressible in $N$.

In this case, the action of
$\pi_1(F)$ on $\cal{T}_F$ is small. By Theorem B (Skora's Theorem),
there are a measured geodesic lamination $\beta$ on $F$ and an
isomorphism $\phi:\cal{T}_\beta\rightarrow \cal{T}_F$ from the dual
tree $\cal{T}_\beta$ of $\beta$ to $\cal{T}_F$. If $\beta$
crosses $\alpha$, then $\alpha$ is realised in
$\cal{T}$ (cf. \cite[Theorem 3.1.4]{meister2}). If $\beta$ and $\alpha$ are
disjoint, then $\alpha$ is collapsed in $\cal{T}$.

It remains to
deal with the case when $\alpha$ is a connected component of
$\beta$. Let $(\eta_i)$ be a sequence of positive numbers converging
to $0$. Let $\kappa_i\subset F$ be a sequence of segments
intersecting $\alpha$ and $\beta$  such that
$\kappa_{i+1}\subset\kappa_i$ and $\int_{\kappa_i}
d\beta\leq\frac{1}{2}\eta_i$. Let $\theta_i$ be a train track
minimally carrying $\alpha$ and having only one switch, which is
$\kappa_i$ (refer to \cite[$\S$ 3.2]{meister1} for the construction
of such a train track). Let $p'$ be a point of $\bigcap_i\kappa_i$,
and $\hat p'\in \bh^2$ a lift of $p'$. This point $\hat p'$
corresponds to a point of $\cal{T}_\beta$ which we shall also denote
by $\hat p'$. Let $p\in \cal{T}$ be the image of $\hat p'$ under
$\phi$. Let $\hat \kappa_i\subset\bh^2$ be the lift of $\kappa_i$
that contains $\hat p'$, and $\hat\theta_i$ the lift of $\theta_i$
that contains $\hat\kappa_i$. Define $\phi_i$ on $\hat\kappa_i$ by
$\phi_i(\hat\kappa_i)=p$. Extend $\phi_i$ to an equivariant map from
the union of the switches of $\hat\theta_i$ to $\cal{T}$. If $\hat b$ is
a branch of $\hat\theta_i$, we define $\phi_i(\hat b)$ to be the
segment of $\cal{T}$ which connects the images of the vertical sides
of $\hat b$. Finally, extend $\phi_i$ to a $\pi_1(S)$-equivariant
map $\phi_i:\bh^2\rightarrow\cal{T}$.

We have thus constructed $\theta_i$ and $\phi_i$, and we need to check that they have the expected properties. Let $\hat b$ be a branch of
$\hat\theta_i$. Translating it by an element of $\pi_1(S)$, we can
assume that $\hat\kappa_i$ contains a vertical side of $\hat b$.
Since $\theta_i$ has only one switch, there is some $g\in\pi_1(S)$
such that $g(\hat\kappa_i)$ contains the other vertical side of
$\hat b$. Let $\hat k_1$ be an arc joining $\hat\kappa_i$ to
$g(\hat\kappa_i)$ whose projection $k_1$ on $S$ lies in $b-|\beta|$.
Then we have $\int_{k_1} d\beta=0$. Let $\hat
k_2\subset\hat\kappa_i$ be an arc joining $\hat p'$ to $\hat k_1$
and let $\hat k_3\subset g(\hat\kappa_i)$ be an arc joining $g(\hat
p')$ to $\hat k_1$. The $\beta$-measures of $k_2$ and $k_3$ are less
than $\int_{\kappa_i} d\beta\leq\frac{1}{2}\eta_i$. Therefore we
have $\int_{k_1\cup k_2\cup k_3} d\beta\leq\eta_i$. This implies
that the distance between $\hat p'$ and $g(\hat p')$ in ${\cal
T}_\beta$ is less than $\eta_i$. It follows then from the
construction of $\phi_i$ that the length of $\phi_i(\hat b)$ is less
than $\eta_i$.

Thus we conclude that when $F$ is incompressible, one
of the conclusions of Lemma \ref{quivabien} holds.\\

\underline{Second case} : $F$ is compressible.

The difference from the previous case is that the action of $\pi_1(F)$ on $\cal{T}_F$ is not small anymore, and hence we cannot use Theorem B directly.
We shall show that if $\alpha$ is neither collapsed nor realised, the following conditions hold: $S(\alpha)$ is incompressible and each component $d$ of $\partial S(\alpha)$ fixes a point in $\cal{T}$ (meaning that each element in the conjugacy class in $\pi_1(N)$ associated to $d$ has a fixed point). These conditions allow us to use Theorem B on $S(\alpha)$ and to construct $\theta_i$ as in the first case.\\

 First we shall use Theorem A to associate a lamination to the action of $\pi_1(F)$ on $\cal{T}$.

Let $(l_n)$ be a sequence of simple closed curves in $S(\alpha)$ converging
to the support of $\alpha$ in the Hausdorff topology.
 By Theorem A (Morgan-Otal's Theorem), there exist, for any $n$, a measured lamination
$\beta_n\in \cal{ML}(S(\alpha))$ and a morphism $\Phi_n:\cal
T_{\beta_n}\ra \cal T_N$, where $\mathcal T_{\beta_n}$ is the $\br$-tree dual to $\beta_n$, such that $\delta_{{\cal
T}_{\beta_n}}(l_n)=\delta_{\cal T_N}(l_n)$.
The difference from the
case when $F$ is incompressible is that $\Phi_n$ may not be an
isomorphism and that $\beta_n$ depends on $l_n$. Let $B$ be the
Hausdorff limit of $|\beta_n|$ after passing to a converging
subsequence. If $\alpha$ is disjoint from $B$, then $\alpha$ is
collapsed in $\cal T_N$. If $\alpha$ intersects $B$ transversely,
then $\alpha$ is realised in ${\cal T}$ as was shown in \cite[Lemma
11]{boches}.
Hence we only need to deal with the case when $|\alpha|$ is contained in $B$.\\

As will be explained below, it follows from results of \cite{boches}, that $B$ can be extended with a homoclinic leaf $h$. Combining this with the assumption that $\alpha$ is a sublamination of a doubly incompressible lamination we shall prove that $S(\alpha)$ is incompressible. Then we will use $h$ to construct an essential $I$-bundle $W\subset N$ whose boundary contains $S(\alpha)$. Using $W$ in the construction of $\beta_n$ will give us enough control on $\beta_n$ to guarantee that $i(\beta_n,\partial S(\alpha))=0$. This will imply that each component of $\partial S(\alpha)$ fixes a point in $\cal{T}$

\begin{cl}      \label{hom}
There is a geodesic lamination $H\supset B$ that contains a homoclinic leaf $h$.
\end{cl}

\begin{proof}
By Claim \ref{limit}, either $\beta_n$ is in tight position with
respect to some meridian in $F$ or there is a homoclinic geodesic
$h_n\subset F$ which does not cross $\beta_n$. By the proof of
\cite[Proposition 2]{boches}, if $\beta_n$ intersects a meridian $m$
and contains no $m$-waves, then $|\beta_n|$ can be extended to a
geodesic lamination with a homoclinic leaf $h_n$. Furthermore by the
proof of \cite[Proposition 1]{boches}, there is a Hausdorff limit of
meridians that does not cross $\beta_n$. (See Appendix, Proposition
\ref{me} and Lemma \ref{noncross}.) Thus we have found in every case
a sequence of meridians in $F$ converging in the Hausdorff topology
to a geodesic lamination $H$ which does not have a transverse
intersection with $B$. By Casson's criterion (Lemma \ref{casson}),
the lamination $H$ contains a homoclinic leaf $h$.
\end{proof}

From now on we assume that $|\alpha|$ is a sublamination of $B$, since we have seen earlier that in the other cases, the conclusion is straightforward. Since $\alpha$ is a sublamination of a doubly incompressible lamination, $h$ can not be a leaf of $\alpha$ (Lemma \ref{anho}).
As we saw in  Lemma \ref{disque}, an incompressible surface cannot
contain a homoclinic geodesic. By construction $F-S(\alpha)$ is incompressible, hence the geodesic
$h$ does not lie in $F- S(\alpha)$.
Since $\alpha$ is arational in $S(\alpha)$, this
implies that there is a half-leaf $h_+$ of $h$ which is asymptotic to a half-leaf of $\alpha$ on $\partial N$. Up to cutting $h_+$, we may assume that $h_+\subset S(\alpha)$. Let $h_-$ be another half-leaf of $h$ which is disjoint from $h_+$ and from $\partial S(\alpha)$.
From the assumption that $\lambda$ is doubly incompressible, we shall deduce that $h_-$ is disjoint from $S(\alpha)$.

\begin{cl}      \label{notasymptotic}
If $|\alpha|\subset B$, then $h_-\subset F-S(\alpha)$.
\end{cl}

\begin{proof}
Seeking a contradiction, we assume that $h_-\subset S(\alpha)$. We shall show that $F-S(\alpha)$ is compressible, which is impossible by construction.

We fix some geometrically finite hyperbolic metric on $N$ and consider the Floyd-Gromov compactification of $\tilde N$.
Let $\tilde h$ be a lift of $h$ to $\tilde N$.
Since $h$ is homoclinic, the two endpoints of $\tilde h$ in $\partial_\infty \tilde N$ coincide.
Let $\tilde h_\pm$ be the lift of $h_\pm$ that lies in $\tilde h$.

Since $\alpha$ is a sublamination of $\lambda$ which is doubly incompressible, then by Lemma \ref{saintptot}, $\tilde h_+$ and $\tilde h_-$ are asymptotic on $\tilde F$. Take a short geodesic arc $\tilde k$ connecting $\tilde h_-$ and
$\tilde h_+$ which does not lie on $\tilde h$ so that they form  a triangle on $\partial \tilde N$ with one vertex at
infinity in $\partial_\infty  \tilde{N}$. Then any lift of a half-leaf of
$\alpha$ entering this triangle must have the same endpoint in $\partial_\infty \tilde N$ as the
lifts of $\tilde h_-, \tilde h_+$ since $h$ does not intersect $\alpha$ transversely.
Therefore such half-leaf of $\alpha$ is trapped
between $\tilde h_-$ and $\tilde h_+$. It follows that we can push $\tilde k$ towards the end of $\tilde h_\pm$ without changing the intersection number
with $\tilde\alpha$, where $\tilde \alpha$ is the preimage of $\alpha$. Thus we obtain  a sequence of geodesic arcs $(\tilde k_i)$, whose lengths tend to zero and whose $\alpha$-measure is $\int_{\tilde k_i}d\tilde\alpha=\int_{\tilde k}d\tilde\alpha$. We project it to $\partial N$ and pass to a subsequence so that $(k_i)$ converges to a point in the Hausdorff topology.
Since the transverse
measure of $\alpha$ is non-atomic, the only way this can happen is that
$\int_{\tilde k} d\tilde \alpha =0$,  i.e., $\tilde \alpha$ lies outside the triangle and $k$ is disjoint from $\alpha$.

Let $\tilde k'$ be the segment on $\tilde h$ between the two endpoints of $\tilde k$ lying on $\tilde h$.
Then $\tilde k$ is homotopic to $\tilde k'$ in $\tilde N$ since $\tilde N$ is simply connected.
Let $k, k'$ be the projections of $\tilde k, \tilde k'$ to $\partial N$, and $m$ be the closed geodesic homotopic to $k \cup k'$.
We see that $k$ is homotopic to $k'$ in $N$ but not on $\partial N$ since both of them are geodesic arcs.
Therefore, $m$ bounds a (possibly singular) disc in $N$, and it is disjoint from $\alpha$ since both $h$ and $k$ are disjoint from $\alpha$.
This contradicts the fact that $F- S(\alpha)$ is incompressible by the Loop Theorem.
\end{proof}

It is quite easy to deduce from this Claim that $S(\alpha)$ is incompressible.

\begin{cl}
If $|\alpha|\subset B$, then $S(\alpha)$ is incompressible.
\end{cl}

\begin{proof}
Let us assume the contrary. By Lemma \ref{tightposition}, there is a meridian $m\subset S(\alpha)$ such that $S(\alpha)$ contains no $m$-waves disjoint from $\alpha$. In particular $h_+$ is in tight position with respect to $m$. Let $(\tilde m_n)\subset\tilde N$, $n\in\bn$, be the family of lifts of $m$ which intersect $\tilde h_+$. Since $h_+$ is in tight position with respect to $m$, each $\tilde m_n$ intersects $\tilde h_+$ once, the curves $\tilde m_n$ are nested and converge to the endpoint of $\tilde h_+$ (\cite[Claim 3.4]{espoir}). Since $\tilde h_+$ and $\tilde h_-$ have the same endpoint,  $\tilde h_-$ intersects $\tilde m_n$ for $n$ large enough. In particular $h_-$ intersects $m\subset S(\alpha)$, contradicting Claim \ref{notasymptotic}.
\end{proof}

We denote by $A$ the minimal lamination contained in the closure of $h_-$. By Claim \ref{notasymptotic}, $S(A)$ is disjoint from $S(\alpha)$.
Since $F-S(\alpha)$ is incompressible, so is $S(A)$.
We shall use Proposition \ref{carac}, to construct an essential $I$-bundle whose boundary contains $S(A)$ and $S(\alpha)$.

\begin{cl}      \label{bundle}
If $|\alpha|\subset B$, then there is an essential $I$-bundle $W\subset N$ whose associated $\partial I$-bundle is $S(\alpha)\cup S(A)$. Moreover, $W$ is homeomorphic to $T\times I$ with $T$ homeomorphic to $S(\alpha)$
\end{cl}

\begin{proof}
Let $\widetilde{S(\alpha)}, \widetilde{S(A)}$ be lifts of
$S(\alpha), S(A)$ containing $\tilde h_+,\tilde h_-$
respectively. Since $S(\alpha)$ and $S(A)$ are incompressible, by
Lemma \ref{disque}, the closures $\overline{\widetilde{S(\alpha)}},\overline{\widetilde{S(A)}}$ of $\widetilde{S(\alpha)}, \widetilde{S(A)}$ in $\overline{\tilde N}$ are discs.

Since $\tilde h_+$ and $\tilde h_-$ have the same endpoint in $\partial_\infty \tilde N$, the two discs $\overline{\widetilde{S(\alpha)}}$ and $\overline{\widetilde{S(A)}}$ intersect in $\partial_\infty \tilde N$.
By Proposition \ref{carac}, there is an essential $I$-bundle $W$ embedded in $(N, S(\alpha)\cup S(A))$, and $h_+$ and $h_-$ are contained in $W \cap (S(\alpha)\cup S(A))$ (after cutting off arcs of finite lengths from $h^+$ and $h^-$). Since $\alpha$ and $A$ lie in the closures of $h_+$ and $h_-$ respectively, the corresponding $\partial I$-bundle contains $S(\alpha)\cup S(A)$. Hence $W\cap\partial N=S(\alpha)\cup S(A)$. This can happen only when $W$ is homeomorphic to $T\times I$ with $T$ homeomorphic to $S(\alpha)$.
\end{proof}

Since $|\alpha|$ lies in the closure of $h_+$ and $A$ lies in the closure of $h_-$, for any leaf $\tilde a$ of the preimage of $\alpha$ in $\widetilde{S(\alpha)}$, there is a leaf $\tilde a'$ of the preimage of $A$ with the same endpoints. It follows that we can isotope the bundle structure so that the projections of $A$ and the support of $\alpha$ along the fibres of $W$ coincide.

Now we can prove that each component of $\partial S(\alpha)$ has a fixed point in ${\cal T}_N$.

\begin{cl}
If $|\alpha|\subset B$, then $\delta_{\cal{T}_N}(\partial S(\alpha))=0$.
\end{cl}

\begin{proof}
 We choose a simple closed curve  $c\subset S(A)$ which is not homotopic to a component of $\partial S(A)$. In the construction of $\beta_n$, we add the conditions that $\delta_{{\cal T}_{\beta_n}}(c)=\delta_{\cal T_N}(c)$ and that $\delta_{{\cal T}_{\beta_n}}(d)=\delta_{\cal{T}_N}(d)$ for any component $d$ of $\partial (T \times\{0\}\cup T \times\{1\})$. Notice that this has no effect on the proof of Claim \ref{hom} and all that follows (in particular Claim \ref{bundle}).

If  $\beta_n$ does not intersect $\partial S(A)$ for some $n$, then $\delta_{{\cal T}_{\beta_n}}(\partial S(A))=0=\delta_{\cal{T}_N}(\partial S(A))$. By Claim \ref{bundle}, $\delta_{\cal{T}_N}(\partial S(\alpha))=\delta_{\cal{T}_N}(\partial S(A))=0$, and we are done.

Otherwise, $B$ intersects $S(A)$. Since $h$ does not cross $B$ and is asymptotic to $A$, this implies that $A$ lies in $B$. In particular, any leaf of $B$ intersecting $\partial S(A)$ contains a half-leaf
asymptotic to $A$. Such a half-leaf intersects
the simple closed curve $c\subset S(A)$ infinitely many times. This
implies that $i(\beta_n,\partial S(A))$ is $o(i(\beta_n,c))$. On
the other hand, by construction, we have
$i(\beta_n,c)=\delta_\cal{T}(c)$ for any $n$. Thus we get
$i(\beta_n,\partial S(A))\longrightarrow 0$.
However, by assumption, $i(\beta_n,\partial S(A))=\delta_{\cal{T}_N}(\partial S(A))$ does not depend on $n$. Thus we have $\delta_{\cal{T}_N}(\partial S(A))=0=\delta_{\cal{T}_N}(\partial S(\alpha))$.
\end{proof}

It follows that the conjugacy class in $\pi_1(F)$ represented by each component of $
\partial S(\alpha)$ has a fixed point in $\cal{T}_N$. This enables us to use Theorem B
on the minimal subtree of $\cal{T}_N$ which is invariant under the action of $\pi_1(S(\alpha))$. Then we can use the same arguments as in the First case to construct the train tracks $\theta_i$ and maps $\phi_i$ and show that the second alternative in the statement of Lemma \ref{quivabien} holds.\\


Thus we have proved the first part of Lemma \ref{quivabien}.
It only
remains to show that at least one component of $\lambda$ is realised
in $\cal{T}$. This was already proved in \cite[Proposition 6.1]{masurext}. Let us briefly review the
proof. Let $(L_n)$ be a sequence of multi-curves converging to
$\lambda$ in the Hausdorff topology. As we have already seen, by Theorem A there
exist a measured geodesic lamination $\beta_n\in\cal{ML}(\partial
M)$ and a morphism $\phi_n:\cal{T}_{\beta_n}\rightarrow \cal{T}$
from the dual tree of $\beta_n$ to $\cal{T}$ such that for any
simple closed curve $l_n$ which lies in $L_n$, either
$\delta_\cal{T}(l_n)>0$ and the restriction of $\phi_n$ to the axis
of $l_n$ is an isometry or $\delta_\cal{T}(l_n)=0$ and
$i(l_n,\beta_n)=0$. Extract a subsequence such that $(\beta_n)$
converges in the Hausdorff topology to a geodesic lamination $B$. As
we have seen above, any connected component of $\lambda$ that
intersects $B$ transversely is realised in $\cal{T}$.

If $S(B)$ is compressible, then, by the proof of \cite[Proposition 2]{boches}, $S(B)$
contains a homoclinic geodesic $h$ which does not cross $B$. Such a homoclinic
leaf crosses $\lambda$ by Lemma \ref{anho}. Thus, if $S(B)$ is compressible, $\lambda$ crosses $B$.

If $S(B)$ is incompressible, then we can apply Skora's Theorem B to
each component of $S(B)$. It follows that $\beta_n$ does not depend
on $n$ for sufficiently large $n$. Denote by $\beta$ this constant
geodesic measured lamination $\beta_n$. We deduce then from
\cite{mors3} that $\beta$ is an annular lamination (see \cite[D\'{e}monstration du Lemme
14]{meister1}). Thus we see that
$\lambda$ crosses the support $B$ of
$\beta$  in this case as well.
\end{proof}

\section{Mapping the train tracks to $\eps_n\bh^3$}\label{track}
In this section, we  shall begin the proof of Theorem \ref{th&}.
Here, we consider a situation more general than the setting in Theorem \ref{th&} as we shall explain below.
We  consider a sequence of measured laminations $(\lambda_n)$ on $S$ converging to $\lambda$.
We also assume that $|\lambda_n|$ converges to a geodesic lamination $L_\infty$ in the Hausdorff topology  and that each minimal sublamination of $L_\infty$ satisfies the conclusion of Lemma \ref{quivabien}.
Let $L_{rec}$ be the union of the recurrent leaves of $L_\infty$.

We shall use the train tracks obtained in lemma \ref{quivabien} to construct a sequence of train tracks $\tau_n$ carrying $\lambda$ and $L_\infty$ which are decomposed into three parts.
The first part $\tau^1$ is independent of $n$ and carries the realised part of $L_{rec}$.
The second part $\tau^2_n$ carries the non-realised part of $L_{rec}$.
The third part $\tau^3_n$, which is $\tau_n -(\tau^1 \cup \tau_n^2)$ and not a train track, carries the rest of $L_\infty$, and has weights with smaller order. Let us denote by $\hat\tau_n$ the preimage of $\tau_n$ under the covering projection from the universal cover of $\partial M$ to $\partial M$.
In the process of constructing $\tau_n$, we shall also build  $\rho_n$-equivariant maps $\tilde f_n$ from $\hat \tau_n$ to $\bh^3$ which map the branches of $\hat\tau^1_n$ to long segments and the branches of $\hat\tau^2_n$ to much shorter ones.  These properties make it possible to estimate the lengths of measured laminations carried by these train tracks.

Recall that when $\hat{\tau}\subset\tilde S$ is the preimage of a
train track $\tau$ on a component $S$ of $\partial M$, we say that a
map $\hat{h}: \hat{\tau} \ra \bh^3$ is $\rho_n$-equivariant if and
only if for every $g \in \pi_1(S)$ and $x \in \hat{\tau}$, we have
$\hat{h}(gx)= \rho_n(i_*(g))\hat{h}(x)$, where $i$ denotes the
inclusion from $S$ to $M$.\\

We shall work in the following setting: $(\rho_n)$ is a sequence of geometrically finite representations uniformising $M$. There is an ultrafilter $\omega$ and $\eps_n\longrightarrow 0$ such that the action of $\rho_n$ on $\eps_n\bh^3$ tends to a small minimal action of $\pi_1(M)$ on an $\br$-tree $\cal{T}$ with respect to $\omega$. We have a sequence $(\lambda_n)$ of measured geodesic laminations converging to a measured lamination $\lambda$.  Furthermore $|\lambda_n|$ converges in the Hausdorff topology to a geodesic lamination $L_\infty$ and each minimal sublamination of $L_\infty$ satisfies the conclusion of Lemma \ref{quivabien}.

We call these assumptions the {\em light assumptions}. We use the assumption that each minimal sublamination of $L_\infty$ satisfies the conclusion of Lemma \ref{quivabien} rather than the one that $\lambda$ is doubly incompressible because it will not contradict the assumption that $l_{\rho_n}(\lambda^*_n)$ is bounded which we shall add in the next section. This way the statements of lemmas and claims in Sections \ref{asum} and \ref{annuli} is non-empty, i.e. not  based on contradictory hypothesis. To illustrate this statement, let us construct an example with these light assumptions which has bounded $l_{\rho_n}(\lambda^*_n)$.

Suppose that $M$ is a handlebody which we regard as an $I$-bundle
$T\times I$ over a compact surface with boundary $T$.  We pick a
convex cocompact representation $\rho$ uniformising $M$, a meridian
$m$, and two pseudo-Anosov diffeomorphisms $\varphi, \psi
:T\rightarrow T$. We extend $\varphi$ and $\psi$ to fibred
diffeomorphisms $\varphi_M,\psi_M:M\to M$. We set
$\rho_n=\psi_{M*}^n\circ\rho$. This sequence tends to the action of
$\pi_1(M)\approx\pi_1(T)$ on the $\br$-tree dual to the stable
lamination of $\psi$.  We set  $m_n=\varphi_M^n(m)$, pick a simple
closed curve $c\subset \partial M$ which is not entirely contained
in $T\times \{0\}$ nor in $T\times \{1\}$, and consider the sequence
$\lambda_n=D_{m_n}^{p_n}\circ\psi_M^n(c)$ where $D_{m_n}:\partial
M\rightarrow\partial M$ is the right Dehn twist along $m_n$. If
$p_n$ is large enough, compared to $n$, the Hausdorff limit
$L_\infty$ of $\lambda_n$ contains two minimal sublaminations both
of which are projected along the fibres to the stable lamination of
$\varphi$. Each of such sublaminations satisfies the conclusion of
Lemma \ref{quivabien} (it is easy to extend each of them to a doubly
incompressible lamination abandoning the other). On the other hand
$\lambda_n$ is homotopic to $\psi_M^n(c)$, hence
$l_{\rho_n}(\lambda^*_n)=l_\rho(c^*)$ is bounded.

We detail the properties of $\tau_n$ and $\tilde f_n$ in the following lemma. We shall work under the light assumptions which were introduced earlier. The rest of the section is devoted to the proof of Lemma \ref{rtt}, which consists in constructing $\tau_n$ and $\tilde f_n$.

\begin{lemma}  \label{rtt}
Under the light assumptions, after taking a subsequence of $(\rho_n)$, there are train
tracks $\tau_n$ carrying both $L_\infty$ and $\lambda_n$ with  decomposition
$\tau_n=\tau^1\cup\tau^2_n\cup \tau_n^3$ such that
\begin{itemize}
\item each minimal sublamination of $L_\infty$ is carried by $\tau^1$ or $\tau^2_n$;
\item $\lambda_n$ is minimally carried by $\tau_n$;
\item $\tau^1$ and $\tau^2_n$ are disjoint subtracks of $\tau_n$
(we do not require $\tau^3_n$ to be a subtrack);
\item $\tau^3_n$ is the union of all branches of $\tau_n$ that are not contained in $\tau^1 \cup \tau^2_n$;
\item the switches of $\tau_n$ lie in $\tau^1\cup\tau^2_n$;
\item the sum of the weights with which $\tau_n$ carries $\lambda_n$ is bounded
independently of $n$;
\end{itemize}
and there is a $\rho_n$-equivariant map $\tilde f_n : \hat\tau_n\rightarrow\bh^3$ from
the preimage of $\tau_n$ in the universal cover of $
\partial M$ to $\bh^3$, which has decomposition $\hat \tau_n =\hat \tau^1 \cup \hat \tau^2_n \cup \hat \tau^3_n$ corresponding to the decomposition $\tau_n=\tau^1\cup\tau^2_n\cup \tau_n^3$,  such that:
\begin{enumerate}[\indent\indent a)]
\item for any branch $\hat b$ of $\hat\tau_n$, its image $\tilde f_n(\hat b)$ is either
a geodesic segment or a point;
\item there are $R>0$ and $n(R) \in \bn$ such that, for $n\geq n(R)$, if
$\hat b$ is a branch
of $\hat\tau^1$, then $l(\tilde f_n(\hat b))\geq R
\eps_n^{-1}$, where $\eps_n$ is the rescaling factor that appeared in the light assumptions;
\item there is a sequence of positive numbers $\delta_n\longrightarrow 0$ such that,
for any $n\in\bn$,
if $\hat b_1,\hat b_2$ are two adjacent branches of $\hat\tau^1$ which are separated by a switch, then the exterior angle between $\tilde f_n(\hat b_1)$ and $\tilde f_n(\hat
b_2)$ is less than $\delta_n$;
\item there is a sequence of positive numbers $\eta_n\longrightarrow 0$ such that, for
any $n\in\bn$, if
$\hat b$ is a branch of $\hat\tau_n^2$, then we have $\eps_n
l(\tilde f_n(b))\leq\eta_n$;
\item for any $n\in\bn$, if $\hat b$ is a branch of $\hat\tau_n^3$, then $\lambda_n(b)(\eps_n l(\tilde f_n(\hat b)))$ is less
than $\eta_n$ for $(\eta_n)$ given in (d).
\end{enumerate}
\end{lemma}

\begin{proof}
It was  proved by Thurston that the set of weighted
simple closed curve is dense in $\cal{ML}(S)$ (cf. \cite{flp} and \cite{penner}), for every component
$S$ of $\partial M$.
By approximating each $(\lambda_n)$ by a
sequence of such unions of weighted simple closed curves, and by a
diagonal extraction, we get a sequence of unions of weighted simple
closed curves satisfying the
light assumption. Therefore we can assume that for
each component $S$ of $\partial M$, the intersection $\lambda_n \cap
S$ is either empty or  a weighted simple closed curve.

Since we have only to construct train tracks on each component of $\partial M$ with
non-empty intersection with $L_\infty$, we can assume that $L_\infty$ is
contained in a component $S$ of $
\partial M$.
 Let $L$ be a minimal sublamination of $L_{rec}$.\\

Let us first consider the case when there is a train track $\theta$ minimally
carrying $L$ which is realised in $\cal{T}_S$ (recall that $\cal{T}_S$ is the minimal subtree of $\cal{T}$ invariant under $i_*\pi_1(S)$).
 Let
$\hat\theta$ be a lift of $\theta$ to the universal cover $\bh^2$ of
$S$. There is a continuous $\pi_1(S)$-equivariant map
$\phi_n:\bh^2\rightarrow\cal{T}_S$ under $\rho_n$ such that $\phi_n$
is constant on every tie of $\hat\theta$ and the restriction of
$\phi_n$ to any rail is injective.
Following \cite{conti}, we
shall construct a $\rho_n$-equivariant map $\tilde
f_n:\hat\theta\rightarrow\bh^3$. Let $\kappa_1,...,\kappa_p$ be the
switches of $\theta$ and $\hat\kappa_1,...,\hat\kappa_p\subset
\bh^2$ lifts of $\kappa_1,...,\kappa_p$.
Denote by $\tilde
x_{i,n}\in\cal{T}_S\subset(X_\omega,x)$ the point
$\phi_n(\hat\kappa_i)$. We first define $\tilde f_n$ on
$\{\hat\kappa_1,...,\hat\kappa_p\}$ by setting $\tilde
f_n(\hat\kappa_i)=\tilde x_{i,n}$. We extend this map to
$\pi_1(S)(\{\hat\kappa_1,...,\hat\kappa_p\})$ by $\tilde
f_n(g(\kappa_i))=\rho_n(g)\circ\tilde f_n(\hat\kappa_i)$ for any
$g\in\pi_1(S)$ and any $1\leq i\leq p$. Let $\hat b$ be a branch of
$\hat \theta$. The vertical sides of $\hat b$ lie in two switches
$\hat\kappa$ and $\hat\kappa'$ whose images by $\tilde f_n$ have
already been defined. On $\hat b$, we let $\tilde f_n$ be the map
which is constant on each tie of $\hat b$, and which induces a
parametrisation of the geodesic segment joining $\tilde f_n(\hat
\kappa)$ to $\tilde f_n(\hat\kappa')$ with constant speed on a
horizontal side of $\hat b$. Then for any branch $\hat b$ of
$\hat\theta$, we have \begin{eqnarray*}\label{limit1} \eps_n l(\tilde
f_n(\hat b))\longrightarrow l_{\cal{T}_S}(\phi_n(\hat b))>0.
\end{eqnarray*}

Let $\theta'$ be the first subdivision of $\theta$ as defined in
\cite[Chapitre 4, \S4.1]{meister2}: that is, for each branch $b_j$ of $\theta$, starting from every endpoint of $b_j \cap b_k$ lying on the vertical side of $b_j$, we cut $b_j$ along a horizontal arc up to its midpoint and reorganise the decomposition into branches keeping the condition that it carries $\lambda_n$.
We let $\hat\theta'\subset\hat\theta$ a lift of $\theta'$.
We shall deform
the map $\tilde f_n$ to one which is adapted to $\hat\theta'$. For a branch $\hat b$ of $\hat\theta'$, its
image by $\tilde f_n$ is a broken geodesic segment which is the union of two geodesic segments. We deform $\tilde f_n$ by
a homotopy which is constant on the vertical sides of $\hat b$ to a
map which is constant on each tie of $\hat b$ and takes $\hat b$ to
the geodesic segment joining the images under $\tilde f_n$ of the
vertical sides of $\hat b$.
By
slightly abusing notation, we shall denote the deformed map by the same
symbol $\tilde f_n$.

Since $\theta$ is realised in $\cal{T}_S$, it follows from the
argument of \cite[Chapitre 4]{meister2} that $\tilde f_n$ has the
following properties:
\begin{enumerate}[\indent\indent a)]
\setcounter{enumi}{1}
\item there are $R>0$ and $n(R)$ such that, for $n\geq n(R)$, if $\hat b$ is a branch
of $\hat\theta'$, we have $l(\tilde f_n(\hat b))\geq R\eps_n^{-1}$;
\item there is a sequence of positive numbers $\delta_n\longrightarrow 0$ such that if
$\hat b_1,\hat
b_2$ are two adjacent branches of $\hat\theta$ which are separated by a switch, then the
external angle between $\tilde f_n(\hat b_1)$ and $\tilde f_n(\hat b_2)$ is smaller than
$\delta_n$.
\end{enumerate}
We repeat the same construction for all the components of $L_{rec}$ that are
realised in $\cal{T}_S$. Denote by $\tau^1$ the union of the train tracks $\theta'$
thus obtained, by $\hat \tau^1$ its preimage in $\bh^2$, and by $\tilde f_n:\hat\tau^1\rightarrow\bh^3$ the map which agrees with the
map defined above on each connected component of $\hat\tau^1$. By Lemma \ref{quivabien},
$\tau^1$ is not empty.
We also see that $\lambda_n$ passes through every branch of $\tau^1$ for every $n$ after taking
a subsequence.\\

When a component $L$ of $L_{rec}$ is not realised in $\cal{T}$, Lemma
\ref{quivabien} gives rise to a sequence of train tracks $\theta_i$ each
carrying
$L$ minimally.
We can assume that $\lambda_i$ passes through every branch of $\theta_i$.
Let us denote by $\tau_i^2$ the union of the train tracks thus obtained from
the
components of $L_{rec}$ which are not realised in $\cal{T}$. Finally we add branches to
$\tau^1\cup\tau^2_i$ to get a train track $\tau_i$ which minimally carries $L_\infty$
and $\lambda_i$.

We shall now extend the map $\tilde f_n$ to the preimage $\hat \tau_i^2$ of $\tau_i^2$. Consider a
connected component $L$ of $L_{rec}$ which is not realised in $\cal{T}$. Consider the
subtrack $\theta_i$ of $\tau_i^2$ which minimally carries $L$. We get from Lemma
\ref{quivabien} that there are a point $x\in\cal{T}$, a sequence $\eta_i\longrightarrow
0$, and
a sequence of $\pi_1(S)$-equivariant maps $\phi_i:\bh^2\rightarrow\cal{T}$ such that
$\phi_i$ maps each branch of the preimage of $\theta_i$ to a geodesic segment (which may be
a point) with length smaller than $\eta_i$ and a lift of $\kappa_i$ (the switch of $\theta_i)$ is mapped to $x$
under $\phi_i$.
Set $x=(\tilde x_n)\in\Pi_n(\eps_n\bh^3)$ (we pick an element in the equivalence class defined by $x$) and fix some $i\in\bn$.
Let $\hat\theta_i\subset\bh^2$ be a lift of $\theta_i$, and $\hat\kappa_i\subset\hat\theta_i$ the lift
of $\kappa_i$ which is mapped to $x$ by $\phi_i$.

Let $G_i\subset\pi_1(\partial M)$ be a finite set consisting of all $g \in \pi_1(\partial M)$ such
that if $\hat b$ is a branch of $\hat\theta_i$ and $\hat\kappa_i$
contains a vertical side of $\hat b$, then the other vertical side
of $\hat b$ lies in $g(\hat\kappa_i)$.
Recall that for each branch
$\hat b$, the length of $\phi_i(\hat b)$ is less than $\eta_i$.
Therefore, we have $d(x,gx)\leq\eta_i$ for any $g\in G_i$. Since
$\cal{T}$ is the $\omega$-limit of $\eps_n \bh^3$, we have $\eps_n
d(\tilde x_n,\rho_n(g)(\tilde x_n))\leq 2\eta_i$ for any $g\in G_i$
for $n$ large enough.
For any $g\in G_i\cup \{id\}$, we define
$\tilde f_n(g(\hat\kappa_i))$ by $\tilde
f_n(g(\hat\kappa_i))=\rho_n(g)(\tilde x_n)$. Let $\hat b$ be a
branch of $\hat\theta_i$ with vertical sides lying in two switches
$\hat\kappa_i$ and $g(\hat\kappa_i)$ for some $g\in G_i$. If $\tilde
f_n(\hat\kappa_i)=\tilde f_n(\rho_n(g)(\hat\kappa_i))$, then we set
$\tilde f_n(\hat b)=\tilde f_n(\hat\kappa_i)$. Otherwise, we set
$\tilde f_n$ to be the map which is constant on each tie of $\hat b$
and takes $\hat b$ to the geodesic segment joining $\tilde f_n(\hat
\kappa_i)$ to $\tilde f_n(\rho_n(g)(\hat\kappa_i))$. Extend $\tilde
f_n$ to an equivariant map from $\hat\theta_i$ to $\bh^3$. For
sufficiently large $n$ and any branch $\hat b$ of $\hat\theta_i$, we
have $\eps_n l(\tilde f_n(b))\leq 2\eta_i$. Furthermore, by construction, the sum of
the weights with which $\theta_i$ carries $\lambda_n$ is bounded by
$\int_{\kappa_i}d\lambda_n\leq\int_{\kappa_1}d\lambda_n\longrightarrow
\int_{\kappa_1} d\lambda$.

We do the same construction for all the components of $L_{rec}$ that are not realised in
$\cal{T}_S$, and denote by $\tilde f_n:\hat\tau^1\cup\hat\tau^2_i\rightarrow\bh^3$ the
maps whose restriction to each connected component of
$\hat\tau^1 \cup \hat\tau^2_i$ is the maps thus defined. It follows from the
construction that there is $N(i)$ such
that for $n\geq N(i)$, if $\tilde b$ is a branch of $\hat \tau^2_i$, we have $\eps_n
l(\tilde f_n(\tilde b))\leq 2\eta_i$.\\

We set $\tau^3_i$ to be the closure of $\tau_i -(\tau^1 \cup \tau^2_i)$, and $\hat \tau^3_i$ its preimage in $\bh^2$.
It remains to define $\tilde f_n$ on  the branches of
$\hat \tau_i^3$. Let $\hat b$ be a branch of $\hat \tau_i^3$. Let $\hat \kappa$ and $\hat
\kappa'$ be the two vertical sides of $\hat b$.
Their projections $\kappa$ and $\kappa'$
lie in $\tau^1\cup\tau^2_i$.
Hence their images by $\tilde f_n$ are already defined, and
there are two points $x=(\tilde x_n), x'=(\tilde x'_n)$ in $\cal{T}$ such that $\tilde
f_n(\bar \kappa)=\tilde x_n$ and $\tilde
f_n(\bar\kappa')=\tilde x'_n$. We set $\tilde f_n$ to be the map which is
constant on each tie of $\hat b$ and takes $\hat b$ to the geodesic segment joining
$\tilde x_n$ to $\tilde x'_n$. We then have $\eps_n d(\tilde
x_n,\tilde x'_n)\longrightarrow d(x,x')$. Furthermore, since $\lambda$
is carried by $\tau^1\cup\tau^2_i$, we have $\lambda_n(b)\longrightarrow 0$.
Therefore, for $n$
large enough, we have $\lambda_n(b)(\eps_n l(\tilde f_n(\tilde b)))\leq 2\eta_i$.

Thus we have proved that there is $N(i)$ such that for $n\geq N(i)$, for a
branch $b$ of $\tau_i-\tau^1$, either $b$ is a branch of $\tau^2_i$ and we have
$\eps_n
l(\tilde f_n(\hat b))\leq 2\eta_i$ or $b$ is a branch of $\tau_i^3$ and we
have $\lambda_n(b)(\eps_n l(\tilde f_n(\hat
b)))\leq 2\eta_i$.
Now by choosing $N(i)$ such that $N(i)< N(i+1)$,
and taking a subsequence $\lambda_{N(n)}$ so that the $n$-th term is
the original $N(n)$-th term, we obtain the desired train track.
This concludes the proof of Lemma \ref{rtt}.
\end{proof}



\section{Finding backtracking} \label{asum}
In this section, we are going to show that for large enough $n$, the path $f_n(c_n)$ has long segments in which it comes back nearly parallel to itself.  Eventually these close returns will allow us to construct some long and thin strips connecting two segments of $f_n(c_n)$. Here we use the adjective \lq\lq long" not only for the rescaled metric of $\eps_n\bh^3$ but also to mean combinatorially long in the sense that they go through many branches of $\tau_n$. Let us fix some notations and definitions to put this idea into a precise statement.

We consider the train tracks $\tau_n$ and maps $\tilde f_n$ which come from Lemma \ref{rtt}.
Take a component $S$ of $\partial M$ such that $S \cap \tau^1 \neq \emptyset$.
In this section, we only have to pay attention to the behaviour of $\lambda_n$ on $S$.
Therefore, for simplicity, we denote $\lambda_n \cap S$ by $\lambda_n$, and $\tau_n \cap
S$ by $\tau_n$, etc. Furthermore, $\lambda_n$ is assumed to be a weighted simple closed curve with support $c_n$ and weight $w_n$.

Let $f_n:\tau_n\longrightarrow\bh^3/\rho_n(\pi_1(M))$ be the projection of
$\tilde f_n$. By construction, $\tau_n$ carries $c_n$. We set
$c^1_n=c_n\cap\tau^1$, $c^2_n=c_n\cap\tau^2_n$ and
$c^3_n=c_n\cap\tau_n^3$. We set $\bar c_n^j=f_n(c_n^j)$ for $j=1,2,3$ and $\bar c_n=f_n(c_n)=\bar c_n^1\cup\bar c_n^2\cup\bar c_n^3$.

Fix an orientation on $c_n$. Let $s$ be a segment lying in $c_n^1$ with the orientation induced by that on $c_n$. The {\em
train route} $b(1),...,b(t)$ {\em of} $s$ is the ordered finite sequence
of branches of $\tau^1$ through which $s$ passes: $b(i)$
is an element of the set $\{b_1,..., b_p\}$ of branches of $\tau^1$.
We fix an orientation for each branch of $\tau^1$. A branch $b(i)$
in the train route of $s$ is said to be {\em positively oriented} if
its orientation coincides with the orientation of $s$ and {\em
negatively oriented} otherwise. The {\em oriented train route}
$bo(1),...,bo(t)$ of $s$ is the ordered finite sequence of oriented
branches of $\tau^1$ through which $s$ goes in this order with the
assigned orientations: $bo(i)$ is an element of $\{b_1,...,
b_p\}\times\{+,-\}$. When $(bo(i))_{i\in I}$ is an oriented train
route, we shall denote by $(b(i))_{i\in I}$ the corresponding
non-oriented train route.

In the following lemma, we shall show that $f_n(c_n)$ nearly
backtracks along some long path. Using the terms of oriented train
routes, this is expressed as follows.

\begin{lemma}
\label{band}
Under the light assumptions, if $(l_{\rho_n}(\lambda_n^*))$ is bounded, then there are two infinite oriented train routes
$bo,bo':\bn\rightarrow \{b_1,..., b_p\} \times \{+, -\}$ in
$\tau^1$ and functions $T,V:\bn\rightarrow\bn$ such that for any $n\in\bn$ there are two disjoint segments $s_n,s'_n\subset c_n^1$ satisfying the following:
\begin{itemize}
\item $T$ is non-decreasing and unbounded;
\item the oriented train routes of $s_n$ and $s'_n$ are $(bo(i))_{0\leq i\leq T(n)}$ and\linebreak $(bo'(i))_{0\leq i\leq V(T(n))}$ respectively;
\item there is a homeomorphism $g_n: f_n(s_n)\rightarrow f_n(s'_n)$;
\item $g_n(f_n(v(i)))\in f_n(b'(V(i)))$ for any $i\leq T(n)$ where
$v(i)=b(i)\cap b(i+1)$;
\item any point $x\in f_n(s_n)$ is connected to $g_n(x)$ by an essential arc
$\zeta_n(x)\subset M_n$ with length less than $6\eps$;
\item the simple closed curve $f_n(s_n)\cup \zeta_{n}(f_n(\partial s_n))\cup f_n(s'_n)$ bounds a disc $D_n$ containing all the arcs $\zeta_n(x)$ for $x\in f_n(s_n)$.

\end{itemize}
\end{lemma}


Recall that $\lambda_n$ is assumed to be a weighted simple closed curve with support $c_n$ and weight $w_n$. We denote by $c_n^*$ the geodesic representative of $c_n$ in $M_n=\bh^3/\rho_n(\pi_1(M))$. Then we set $l_{\rho_n}(\lambda_n^*)=w_n l_{\rho_n}(c_n^*)$ by definition.

Notice that we have adopted the conventions that $0\in\bn$, and that if $T(n)=0$, there are no segments $s_n$ and $s'_n$.

We shall see that $V(i)\leq i\frac{R'}{R}+1$ for any $i\in\bn$ (equation (\ref{large})).

\begin{proof}

First we shall show that most points $x$ in $\bar c_n^1$ are close to another component of $\bar c_n^1$  (i.e. not the component containing $x$). The proof goes roughly as follows: we construct a simplicial annulus $A_n$ between $\bar c_n$ and $c_n^*$.  If a point $x$ of $\bar c_n^1$ is not close in $A_n$ to a point in another component of $\bar c_n^1$ then either $x$ is close to a component of $\bar c_n^2$ or $\bar c_n^3$ or $x$ is not close to any component of $\bar c_n^j$ for $j=1,2,3$. Using the Gauss-Bonnet formula and the length comparison between the components of $\bar c_n^j$, we shall show that this can happen only for the minority of the points of $\bar c_n^j$.\\

Let us start the formal proof.
The curve $\bar c_n=f_n(c_n)$ is a piecewise geodesic.
We define the edges of $\bar c_n$ to be the images of the intersections of $c_n$ with the branches of $\tau_n$, and
the vertices to be the images of the
intersections of $c_n$ with the switches of $\tau_n$.

Let $x_{n,1},...,x_{n,p_n}$ be the vertices of $\bar c_n$, and choose
the same number of points $y_{n,1},...,y_{n,p_n}$ on $c^*_n$.
We shall make a strip
with boundaries $(x_{n,i})$ and $(y_{n,i})$ and  triangulate it
by making each rectangle into a pair of triangles.
To be more precise, for
$1\leq i\leq p_n$, we consider the geodesic triangle with vertices
$y_{n,i},x_{n,i},x_{n,i+1}$ (with $x_{n,p_{n}+1}=x_{n,1}$ and
$y_{n,p_{n}+1}=y_{n,1}$) and the geodesic triangle with vertices
$x_{n,i+1},y_{n,i},y_{n,i+1}$.
The union of these triangles for
$i=1, \dots , p_n$ is a simplicial annulus $A_n=S^1\times [0,1]$
bounded by $c^*_n$ and $\bar c_n$.
The metric $\nu_n$ induced on
this annulus by the lengths of paths is a hyperbolic metric with
piecewise geodesic boundary.
By the Gauss-Bonnet formula, the area
of $A_n$ is less than $2p_n\pi$.
By Lemma \ref{rtt}, the sequence
$(w_np_n)=(\sum_{b\text{: the branches of }\tau_n} \lambda_n(b))$ is
bounded. We parametrise $A_n$ by $S^1\times[0,1]$ so that the
projection of $S^1\times\{1\}$ to $\bh^3/\rho_n(\pi_1(M))$ is
$c_n^*$.

For a positive number $\epsilon$, which we shall specify later, and each point
$x\in\bar c_n^1$, we consider a geodesic arc $a_x$ on $(S^1 \times I, \nu_n)$
perpendicular to $S^1 \times \{0\}$ at $x$ having length $\epsilon$ with respect to
$\nu_n$. If the perpendicular reaches $ S^1 \times \partial I$ before the length
$\epsilon$ is attained, we define $a_x$ to be the geodesic arc having both endpoints on $
S^1 \times \partial I$.\\

We shall show that most arcs $a_x$ issuing from a component of $\bar c_n^1$ intersect an arc $a_y$ issuing from another component of $\bar c_n^1$. For this purpose, we are going to estimate from below the length of the set of points
$x$ for which the $a_x$ reach $S^1 \times \{1\}$ without
intersecting themselves or each other. Since the length of this set
of points is bounded above by $\mathrm{length}(c_n^*)$ (with respect
to $\nu_n)$, we get an inequality, which will appear as the
inequality (i) below. For that, we need to subtract from the length
of $\bar c_n^1$ the lengths of (I) the set of points $x$ for which
$a_x$ has self-intersection, (II) the set of points $x$ for which
$a_x$ intersects $a_y$ with $x \neq y$, (III) the set of points $x$
for which $a_x$ has an endpoint on either $\bar c_n^2$ or $\bar
c_n^3$, and (IV) the set of points $x$ which are neither of type (I)
nor of type (II) and for which $a_x$ has an endpoint in the interior
of $S^1\times I$. See Figure \ref{types}.

\begin{figure}[hbtp]
\psfrag{1}{$\bar c_n^1$}
\psfrag{2}{$\bar c_n^2$}
\psfrag{3}{$\bar c_n^3$}
\psfrag{x}{$x_{n,1}$}
\psfrag{v}{$x_{n,2}$}
\psfrag{y}{$y_{n,1}$}
\psfrag{w}{$y_{n,2}$}
\psfrag{c}{$c_n^*$}
\includegraphics[width=.8\linewidth]{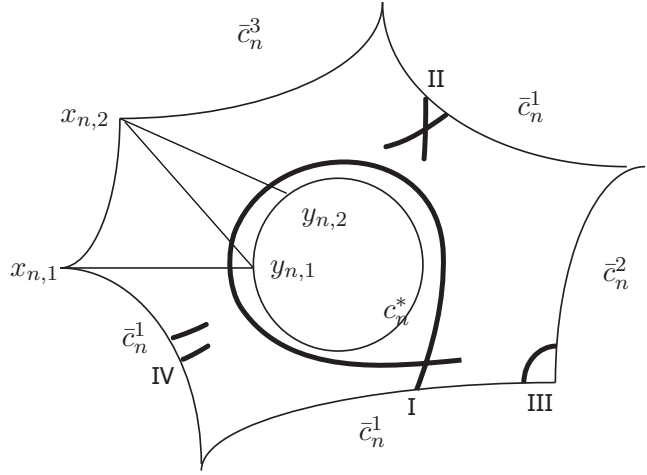}
\caption{Annulus between $c_n^*$ and $f_n(c_n)$}
\label{types}
\end{figure}

We first consider the contribution of the points of type (I) to the length,
i.e., $x$ for which
$a_x$ intersects itself transversely. By the Gauss-Bonnet formula, a geodesic loop
formed by a
subarc of $a_x$ cannot be null-homotopic. Hence, there must be a loop formed by a subarc
of $a_x$ freely homotopic to $S^1 \times \{1\}$. It follows that if both of two
perpendiculars $a_{x_1}, a_{x_2}$ with $x_1 \neq x_2$ have self-intersection, then
$a_{x_1} \cap a_{x_2} \neq \emptyset$. Thus, the contribution of
the set of
$x$ with self-intersecting $a_x$ (i.e. of type (I)) to the length is absorbed in the contribution of $x$ with $a_x$
intersecting another $a_y$, that is, of type (II), which will be dealt with below.

We next consider the points $x$ of type (II),  i.e.\ those for which $a_x$ intersects
$a_y$ for some $y\in\bar c_n^1$.
 Let $m$ be a point in the
intersection $a_x \cap a_y$, and let $a_x', a_y'$ be subarcs of
$a_x, a_y$ between $x$ and $m$ and $y$ and $m$ respectively. Let
$\beta$ be an arc on $\bar c_n$ to which $a_x' \cup a_y'$ is
homotopic fixing the endpoints. Suppose that $\beta$ is also
contained in $\bar c_n^1$. We then say that $x$ is an {\em
inessential} point of type (II) and that $a'_x\cup a'_y$ is an {\em
inessential arc}. It was shown in \cite[Lemme 5.11]{Bouts} that, in
this situation, there is a constant $\xi_n$ depending only on
$\epsilon$ and the maximal exterior angle of the vertices on $\bar
c_n^1$, which is less than $\delta_n$ in our case, such that $x$ is
within distance $\xi_n$ with respect to $\nu_n$ from a vertex of
$\bar c_n^1$. It was also shown in \cite{Bouts} that the constant
$\xi_n$ tends to $0$ as either $\delta_n$ or $\epsilon$ goes to $0$.
If $\beta$ does not lie on $\bar c_n^1$, then we say that $x$ is an
{\em essential} point of type (II). See Figure \ref{type 2}.
\begin{figure}[hbtp]
\psfrag{1}{$\bar c_n^1$}
\psfrag{2}{$\bar c_n^2$}
\includegraphics[width=.8\linewidth]{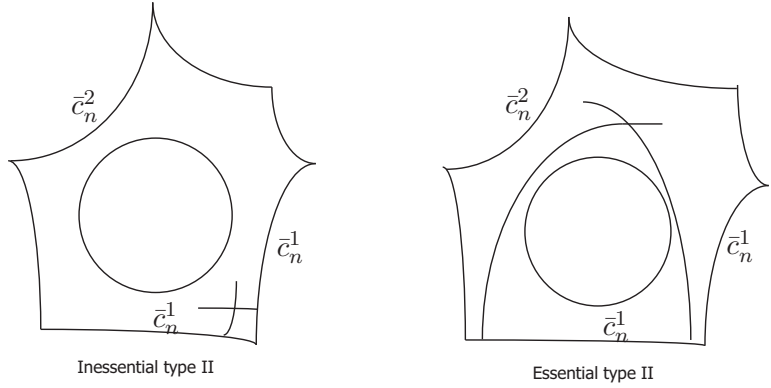}
\caption{Inessential and essential points of type II.}
\label{type 2}
\end{figure}

We note that the
length of an essential arc (such as $a'_x\cup a'_y$)  is less than or equal to
$2\epsilon$. Let $\bar c_n^+$ be the union of the essential points
of type (II). We shall use the essential points of type (II) to
construct the long and thin strips of Lemma \ref{band}.\\

Now we shall bound the length of the sets of points of type (III) and (IV). First we consider the points of type (III). The total length with respect to $\nu_n$
of the set of points $x$ on $\bar c_n^1$ for which $a_x$ reaches a point on $\bar c_n^2$
is bounded above by the length of $\bar c_n^2$. Similarly, the total length of the set of
points $x$ for which $a_x$ reaches a point on $\bar c_n^3$ is bounded above by
the length of $\bar c_n^3$.

Finally we consider the points of type (IV); the points $x$ such that $a_x
\setminus \{x\}$ is contained in the interior of $S^1\times I$ while $a_x$ has
neither
self-intersection nor intersection with another $a_y$. Since the union of $a_x$ for $x$
of type (IV) has area bounded below by the length of the set of points $x$ of type (IV)
multiplied by $\sh(\epsilon)$, we can bound the length from above by
$\mathrm{Area}(A_n)/\sh(\epsilon)$.

Putting all of these considerations together, we get an inequality:
\begin{eqnarray*}
\text{(i)} \hspace{3mm} \l(\bar{c}_n^1) -2p_n\xi_n -\l(\bar{c}_n^+) - \l(\bar{c}_n^2)
-\l(\bar{c}_n^3) -\mathrm{Area}(A_n)/\sh(\epsilon)
\leq \l(c_n^*).
\end{eqnarray*}

\smallskip
Notice that we have
$$w_nl_{\rho_n}(\bar c_n^2)=\sum_{b\text{: branches of
}\tau_n^2}\lambda_n(b)l_{\rho_n}(f_n(b))$$ and
$$w_nl_{\rho_n}(\bar c_n^3)=\sum_{b\text{: branches of
}\tau_n^3} \lambda_n(b)l_{\rho_n}(f_n(b)).$$
Therefore, by the property (d) of Lemma \ref{rtt}, we have
$w_nl_{\rho_n}(\bar c_n^2)=o(\eps_n^{-1})$ and by the property (e),
we have $w_nl_{\rho_n}(\bar c_n^3)=o(\eps_n^{-1})$.

By Lemma \ref{rtt} $(w_np_n)$ is a bounded sequence. It follows that we have $2w_n
p_n\xi_n\longrightarrow 0$. This implies also that
$w_n\mathrm{Area}(A_n)\leq 2 w_n p_n\pi$ is a bounded sequence and
that we have $\epsilon_n
w_n\mathrm{Area}(A_n)\longrightarrow 0$.

By assumption $(w_n
l_{\rho_n}(c^*_n))=(l_n(\lambda_n))$ is a bounded sequence; hence
$\epsilon_n w_n l_{\rho_n}(c^*_n)$ tends to $0$. Thus we have shown
the following.
\begin{cl} \label{zero} We have
$w_n\epsilon_n(l_{\nu_n}(\bar{c}_n^1)-l_{\nu_n}(\bar{c}_n^+))\longrightarrow 0$.\hfill $\Box$
\end{cl}

Now that we know that $\bar c_n^+$ occupies the most part of $\bar c_n^1$, we shall use $\bar c_n^+$ to construct maps $g_n$ and strips $D_n$.
First we define a discrete version of $g_n$.
Let $\{\sigma_{1,n},\sigma_{2,n},\dots \}\subset\bar c_n$ be a maximal family of
disjoint
segments with diameter $6\eps$ such that the midpoint $x_{i,n}$ of $\sigma_{i,n}$ lies
in $\bar c_n^+$. Consider a segment $\sigma_{i,n}$ and its middle point $x_{i,n}$. In the
family of essential arcs joining $x_{i,n}$ to $\bar c_n^1$, we take an arc
$a_{i,n}$
to be the shortest (with respect to $\nu_n$). Since $x_{i,n}$ lies in
$\bar c_n^+$, the length of $a_{i,n}$ is less than $2\eps$.

 If $
\partial a_{i,n}-x_{i,n}$ lies in some $\sigma_{j,n}$, we denote by $\zeta_{i,n}$
the
geodesic segment in $(S^1 \times I, \nu_i)$ joining $x_{i,n}$ to $x_{j,n}$ which is homotopic to $a_{i,n}$ relative
to $x_{i,n}\cup\sigma_{j,n}$. If $
\partial a_{i,n}-x_{i,n}$ is disjoint from $\bigcup\sigma_{j,n}$, we define
$\zeta_{i,n}$
to be $a_{i,n}$ (see Figure \ref{dis}). The length of each segment $\zeta_{i,n}$ thus obtained is less
than $5\eps$.
Using the minimality of the lengths of the arcs $a_{i,n}$, we shall next show that the segments $\zeta_{i,n}$ have mutually disjoint interiors.

Consider two different segments $\zeta_{1,n}$ and $\zeta_{2,n}$ and assume that
their interiors intersect. Then the interiors of $a_{1,n}$ and $a_{2,n}$
also intersect.
Let  $y$ be a point in the intersection. Let $[x_{\ell,n},y[,\ \ell=1,2$ be the connected
component of $a_{\ell,n}-\{y\}$ containing
$x_{\ell,n}$.
Let $\gamma_n$ be the shortest of the two segments $a_{1,n}-[x_{1,n},y[$ and
$a_{2,n}-[x_{2,n},y[$ (see Figure \ref{dis}). Then, for $\ell=1,2$, the length of the arc
$[x_{\ell,n},y]\cup\gamma_n$ is less
than or equal to the length of $a_{\ell,n}$.
Furthermore one of the two arcs
$[x_{\ell,n},y]\cup\gamma_n$, say $[x_{1,n},y]\cup\gamma_n$ is not the shortest in its
homotopy class relative to the endpoints.
Let $a'_1$ be the shortest arc homotopic to $[x_{1,n},y]\cup\gamma_n$ relative
to the endpoints.
Then
the length of the segment $a'_1$ (with respect to $\nu_n$) is less than the length of $a_{1,n}$.
Recall that we
chose $a_{1,n}$ which is shortest among all essential arcs joining
$x_{1,n}$ to $\bar c_n^1$. It follows that $a'_1$ is not essential, i.e.\ there is a
segment $\beta\subset \bar c_n^1$ homotopic to $a'_1$ relative to the endpoints.
The endpoints of $\beta$ are
$x_{1,n}$ and another point which we call $y_1$. The distance (with respect to $\nu_n$)
between $x_{1,n}$ and $y_1$ is less than the length of $a'_1$ which is less than $2\eps$.
By the properties (b) and (c), each component of $\bar c_n^1$ (in particular
the one
containing $\beta$) is a union of long geodesic segments such that the exterior angle
between two consecutive segments is small. By \cite[Lemma 4.2.10]{ceg} such a component of
$\bar c_n^1$ is a $(K, \eta)$-quasi-geodesic with $K \ra 1, \eta \ra 0$ as $n
\ra \infty$.
It follows that there is $N$ (independent of $\beta$)
such that for $n\geq N$, the length of $\beta$ is less than $3\eps$. This implies that
$y_1$ lies in $\sigma_{1,n}$. By our definition of $\zeta_{1,n}$, it has an
endpoint on $x_{1,n}$, not on $y_1$ (see Figure
\ref{dis}). This contradicts our assumption that the interiors of $\zeta_{1,n}$ and $\zeta_{2,n}$ intersect.

\begin{figure}[hbtp]
\psfrag{a}{$x_{1,n}$}
\psfrag{b}{$y_1$}
\psfrag{c}{$y$}
\psfrag{d}{$x_{2,n}$}
\psfrag{f}{$\bar c_n$}
\psfrag{g}{$\sigma_{1,n}$}
\psfrag{h}{$\sigma_{2,n}$}
\psfrag{i}{$c^*_n$}
\psfrag{j}{$a_{2,n}$}
\psfrag{k}{$a_{1,n}$}
\psfrag{l}{$\zeta_{2,n}$}
\psfrag{m}{$\zeta_{1,n}$}
\psfrag{n}{$\gamma_n$}
\centerline{\includegraphics{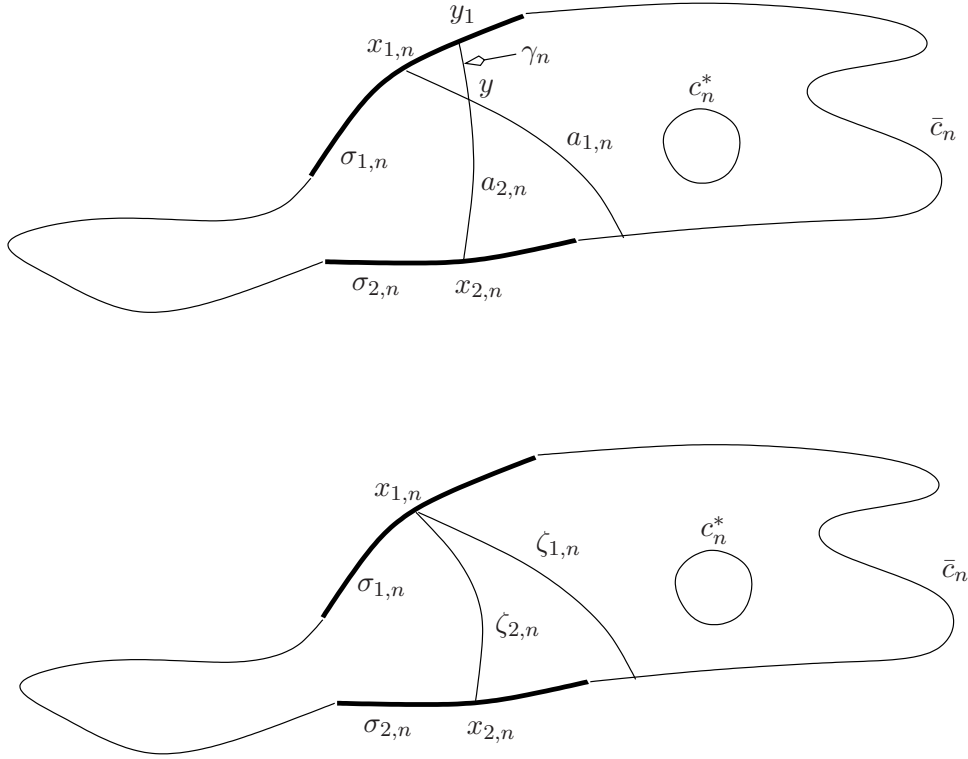}}
\caption{From $a_{i,n}$ to $\zeta_{i,n}$}
\label{dis}
\end{figure}

\indent
Even if some segment $a_{i,n}$ has a self-intersection, the same argument shows
that $\zeta_{i,n}$ does not have any self-intersection.

Thus we have proved the following claim.

\begin{cl} There are a map $h_n:\{x_{1,n},x_{2,n},...\}\rightarrow \bar c^1_n$ and
a family $(\zeta_{i,n})$ of essential segments with disjoint interiors such that the
length of $\zeta_{i,n}$ is less than $5\eps$,  $\partial\zeta_{i,n}=\{x_{i,n},h_n(x_{i,n})\}$, and
$\{h_n(x_{1,n}),h_n(x_{2,n}),...\}\cap\bigcup_i \sigma_{i, n
}\subset\{x_{1,n},x_{2,n},...\}$.\hfill $\Box$
\end{cl}

We shall use those segments $\zeta_{i,n}$ as a skeleton to construct our strips $D_n$. To this end we need to fill the space between two consecutive segments $\zeta_{i,n}$ and $\zeta_{i+1,n}$ with a thin strip. In particular, we need to check that  $\zeta_{i,n}$ and $\zeta_{i+1,n}$ are in the same homotopy class and that each point of $\bar c_n^1$ between $x_{i,n}$ and $x_{i+1,n}$ is close to a point on $\bar c_n$ between $h_n(x_{i,n})$ and $h_n(x_{i+1,n})$.
Our proof will rely mainly on counting segments. Namely we start with sufficiently many segments and prove that at each step,there are only few among them which fail to satisfy the property that we need.\\

Let $j$ be a positive integer, and cut $\bar c^1_n$ into disjoint segments
$\bar  s_{1,n},\bar  s_{2,n},...$, each containing $j$ edges (if the number of edges of some
component of $\bar c^1_n$ is not a multiple of $j$, then there are some edges of
$\bar
c^1_n$ not belonging to any one of these segments).

We shall evaluate the number of segments thus obtained using the following claim.
\begin{cl}
\label{order of number}
Let $r_n$ be the number of the components of $\bar c_n^1$. Then, $w_nr_n
\longrightarrow 0$ as $n \rightarrow \infty$.
\end{cl}

\begin{proof}
 Note that any train route on $\tau_n$ connecting a point in $\tau^1$ and a
point in $\tau^2_n$ must pass through a point in $\tau_n^3$.
Therefore between any two distinct components of $c_n^1$,
there is a component of $c_n^3$.
Hence $(w_nr_n)$ is bounded above by $\sum_{b \subset \tau_n^3}\lambda_n(b)$. Since $|\lambda|\subset
L_{rec}$ is carried by $\tau^1\cup\tau^2_n$, the sum $\sum_{b
\subset \tau_n^3}\lambda_n(b)$ tends to $0$ as
$n \rightarrow \infty$. It follows that we have $w_nr_n
\longrightarrow 0$.
\end{proof}

Since \begin{eqnarray}\label{edge}
w_n (\text{number of edges of}\ \bar c^1_n)=w_n(\sum_{b\subset \tau^1}
c_n(b))
\longrightarrow \sum_{b\subset\tau^1}\lambda(b),
\end{eqnarray} the number of edges of $\bar
c^1_n$ is $\Theta(w_n^{-1})=\Theta(p_n)$ (see Section \ref{notations} for notations), where $p_n$ was defined to be the number vertices on $\bar c_n$ (i.e. the number of times $c_n$ crosses a switch of $\tau_n$). The number of edges of $\bar c_n^1$  lying in none of the $\bar s_{i,n}$ is less than $j r_n=o(p_n)$.
It follows that the number of the segments $\bar s_{i,n}$ is $\Theta(p_n)$.

Let $t_n$ be the number of edges of $\bar c^1_n$ containing no one
among the segments $\sigma_{i,n}$ defined earlier. If an edge $e$
contains no
segment among the $\sigma_{i,n}$, then there is no point of $\bar c_n^+$ in $e$
outside the $3\eps$-neighbourhood of $
\partial e$. By the property (b) in Lemma \ref{rtt}, the total length of
these edges is greater than $t_n R\epsilon_n^{-1}$ and smaller than $l_{\nu_n}(\bar c_n^1)-l_{\nu_n}(\bar c_n^+)+6t_n\eps$. By Claim \ref{zero} and the equation (\ref{edge}), we have $w_n\eps_n(l_{\nu_n}(\bar
c_n^1)-l_{\nu_n}(\bar c_n^+)+6t_n\eps)\longrightarrow 0$. Hence $w_n\epsilon_n t_n
R\epsilon_n^{-1}\longrightarrow 0$ and $t_n$ is  $o(w_n^{-1})=o(p_n)$.

Thus we know that among the $\bar s_{i,n}$, there are $\Theta(p_n)$ disjoint
segments lying in
$\bar c^1_n$ and containing $j$ edges (where $j$ is the number we have fixed when cutting $\bar c_n^1$) each of which contains a segment
among the $\sigma_{i,n}$.
We shall denote these $\Theta(p_n)$-many segments again by $\bar s_{1,n},
\bar s_{2,n}, \dots $.\\

Next we shall show that most of these segments are joined to only one component of $\bar c_n^1$ through the arcs $\zeta_{i,n}$.

Let $\bar s\in\{\bar s_{1,n},\bar s_{2,n},\dots \}$ be a segment with the following
property: there are at least two distinct components of
$\bar c^1_n$ containing an endpoint of $\zeta_{i,n}$ for some $x_{i,n}\in \bar s$.
Let $t'_n$ be the
number of those with this property among the $\bar s_{i,n}$. In each such   segment
$\bar s$, we choose
two points in $\{x_{1,n},x_{2,n},...\}\cap \bar s$, say $x_{1,n},x_{2,n}$, such
that $\zeta_{1,n}$ and $\zeta_{2,n}$ connect $\bar s$ to distinct components of
$\bar c^1_n$ and
the segment $]x_{1,n},x_{2,n}[\subset \bar s$ contains no $x_{k,n}$.
We note that there may be  $\zeta_{i,n}$ other than $\zeta_{1,n}, \zeta_{2,n}$ which have $x_{1,n}$ or $x_{2,n}$ as an endpoint.
There are two points $y_{1,n}$ and $y_{2,n}$ which lie in two distinct
components of $\bar c_n^1$ such that $y_{1,n}$ (resp. $y_{2,n}$) is connected to
$x_{1,n}$ or $x_{2,n}$ by some $\zeta_{i,n}$ and that $y_{1,n}$ and $y_{2,n}$ are innermost in the following sense : if $[y_{1,n},y_{2,n}]$ is the
segment of $\bar c_n$ joining $y_{1,n}$ to $y_{2,n}$ whose interior does not
contain $x_{1,n}$, then there is no $\zeta_{i,n}$ connecting $]y_{1,n},y_{2,n}[$
to $\{x_{1,n},x_{2,n}\}$ (these $y_{1,n}, y_{2,n}$ may or may not coincide with $h_n(x_{1,n}), h_n(x_{2,n})$).
Let us embed $A_n$ into a round disc in such a way that $\bar
c_n$ is the boundary of the disc, and connect $y_{1,n}$ to $y_{2,n}$ by a geodesic segment
with respect to the ordinary Euclidean metric on the disc.

Assume that for another segment $\bar s' \in \{\bar s_{1,n}, \dots\}$ with the same
property, the resulting geodesic segment in the round disc intersects
transversely the geodesic segment produced from $\bar s$ above (i.e. the geodesic segment
connecting $y_{1,n}$ to $y_{2,n}$). Suppose that $x_{3,n}$ and $x_{4,n}$ are the
points on $\bar s'$ chosen in the same way as $x_{1,n}, x_{2,n}$ for $\bar s$.
We number them so that the order in which the four points lie on the circle is $x_{1,n}, x_{2,n}, x_{3,n}, x_{4,n}$.
Since $\zeta_{i,n}$ have disjoint interiors and $y_{1,n}$ and $y_{2,n}$ are innermost, we see that the only ways this can happen are the following two:
(1) $x_{3,n}=y_{1,n}$ and (2) $x_{4,n}=y_{2,n}$.
Therefore for each $\bar s$, there are only two configurations of $\bar s'$ such that the geodesics in the round disc intersect transversely.
See Figure \ref{poly}.
%

\begin{figure}[hbtp]
\psfrag{a}{$x_{2,n}$}
\psfrag{b}{$x_{1,n}$}
\psfrag{c}{$y_{1,n}=x_{3,n}$}
\psfrag{d}{$x_{4,n}$}
\psfrag{e}{$y_{4,n}$}
\psfrag{f}{$y_{2,n}$}
\psfrag{g}{$\bar s$}
\psfrag{h}{$\bar s'$}
\psfrag{i}{$y_{2,n}=x_{4,n}$}
\psfrag{j}{$y_{1,n}$}
\psfrag{k}{$x_{3,n}$}
\psfrag{l}{$\bar c_n$}
\psfrag{m}{$c^*_n$}
\psfrag{o}{Case (1)}
\psfrag{p}{Case (2)}
\centerline{\includegraphics{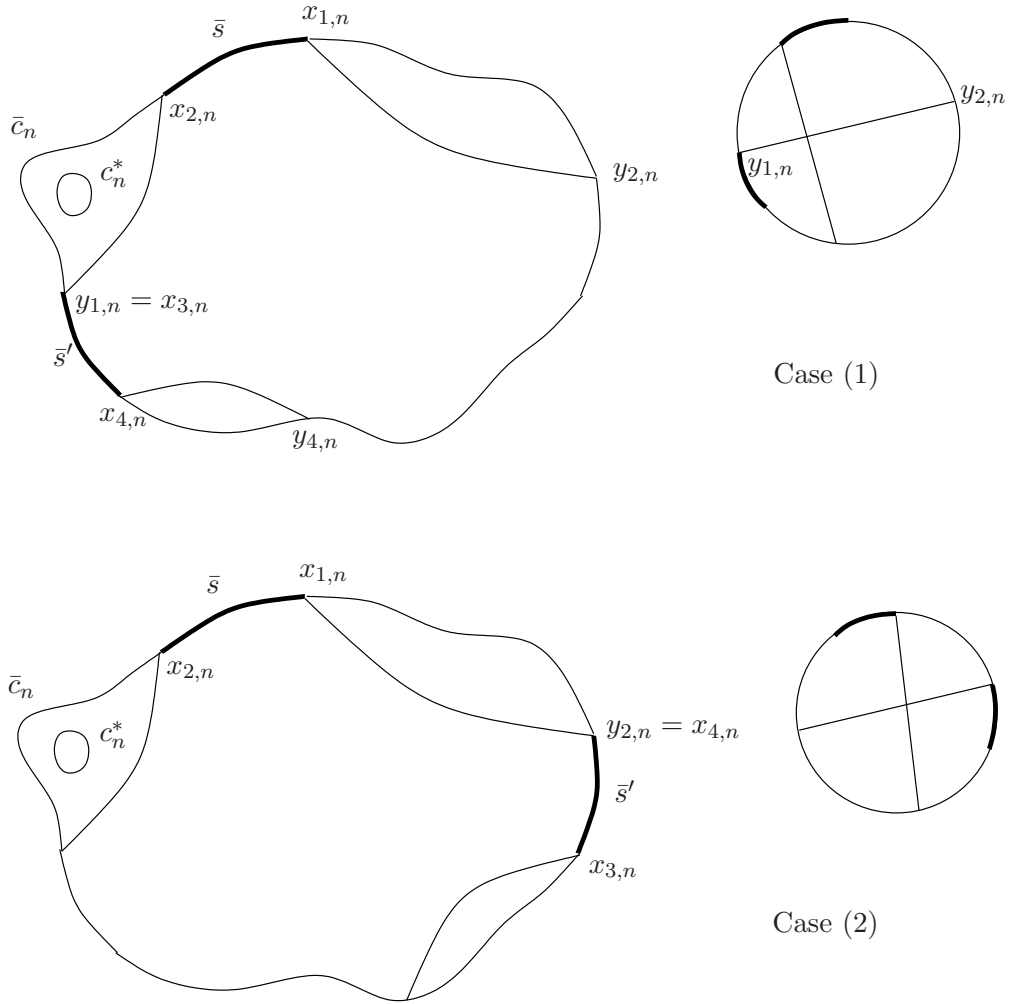}}
\caption{A segment intersects at most two other segments}
\label{poly}
\end{figure}

Thus we can build at least $\frac{1}{3} t'_n$ disjoint
geodesic segments in the round disc, each connecting two distinct components
of $\bar c_n^1$.
Furthermore, since the $\zeta_{i,n}$ have disjoint interiors, any pair of connected
components of $\bar c_n^1$ is connected by at most two of these disjoint
segments.
Consider a map from $\bar c_n$ to a round circle which preserves the order and maps
each connected component $C_i$ of $\bar c_n^1$ to a point $Q_i$. Join two points
$Q_i$ and $Q_j$ by a segment if and only if there is one of the segments constructed above which
joins $C_i$ and $C_j$. Thus we have constructed at least $\frac{1}{6} t'_n$
segments with
disjoint interiors. We can now add some geodesic segments in the disc to get a
triangulation of
the polygon with vertices $Q_1,...,Q_{r_n}$. Such a triangulation has $2r_n-3$
edges (this can easily be computed with the Euler formula). Therefore we have
$\frac{1}{6} t'_n\leq 2r_n-3=o(p_n)$.

Since we initially had $\Theta(p_n)$ segments $\{\bar s_{1,n},\bar s_{2,n},...\}$,
after excluding $o(p_n)$ segments as above from them,
there remains, by abusing notations again, $\Theta(p_n)$ disjoint segments $\{\bar s_{1,n},
\bar s_{2,n},...\}$ in $\bar c_n^1$,  such that
if $\bar s$ is one of those segments, then we have:
\begin{enumerate}[- ]
\item $\bar s$ contains $j$ edges;
\item i) each edge of $\bar s$ contains some $x_{i,n}$;
\item ii) there is a unique component $C$ (depending on $\bar s$) of $\bar c_n^1$ such
that for any $x_{i,n}\in \bar s$, we have $h_n(x_{i,n})\in C$.
\end{enumerate}

Let $\bar s$ be one of these segments, $C$ the associated component of
$\bar c_n^1$, and $x$ a point of $\bar s\cap\{x_{1,n},...\}$.
Denote by
$\zeta_x$ the corresponding segment $\zeta_{i,n}$. Since $A_n$ is an annulus and
$\zeta_x$ is embedded, there are only two possibilities for the homotopy class
of $\zeta_x$ relative to $\bar s\cup C$. Therefore, taking $2j$ instead of $j$ at the
beginning and cutting each segment into two groups, we get
$\Theta(p_n)$ disjoint segments $\{\bar s_{1,n},...\}$ in $\bar c_n^1$, each one
containing $j$ edges and satisfying (i), (ii) above and:\smallskip\\
\indent iii) for any $x,y\in \bar s\cap \bar c_n^+$, $\zeta_{x}$ and $\zeta_{y}$ are
homotopic relative to $\bar s\cup C$.\medskip\\

Now we have a sufficiently number of $j$ consecutive segments $\sigma_{i,n}$ with the properties (i), (ii) and (iii). Next we shall show that they lie in the boundary of the expected strips $D_n$.

\indent
Let $\bar s$ be one of the segments produced above, and $C$ the
corresponding component of $\bar c_n^1$ in the property (ii). Let $x$ and $y$ be the extremal points
of $\bar s\cap\{x_{1,n},...\}$, and  $[x,y]$ the segment in $\bar s$ joining $x$ to $y$.
The segment $[x,y]$ contains at least $(j-2)$ edges. We have the following:

\begin{lemma}  \label{anne} There is $N\in\bn$ which does not depend on $\bar s$ such
that for any $n\geq N$ we have the following by homotoping $A_n$ keeping
$\partial A_n$ and $\bar c_n^i$ unchanged:\\
There is a homeomorphism $g_n:[x,y]\rightarrow [h_n(x),h_n(y)]$ such
that for any $n\geq N$ and any $z\in[x,y]$, the two points $z$ and $g_n(z)$ are
connected by an
essential arc whose length (with respect to the induced metric $\nu_n$ on $A_n$) is less than $6\eps$.
\end{lemma}

\begin{proof} By the property (iii), the simple closed curve $\zeta_x\cup
[x,y]\cup \zeta_y\cup [h_n(x),h_n(y)]$ bounds a disc in $A_n$. Since
both $[x,y]$ and $[h_n(x),h_n(y)]$ lie in $\bar c_n^1$, by
the properties (b) and (c) in Lemma \ref{rtt}, they consist of long geodesic segments such
that the external angles formed by two adjacent segments are less
than $\delta_n$. Let $k$ be the geodesic segment in
$\bh^3/\rho_n(\pi_1(M))$ joining $x$ to $y$ which is homotopic to
$[x,y]$, and let $k'$ be the one in $\bh^3/\rho_n(\pi_1(M))$ joining
$h_n(x)$ to $h_n(y)$ which is homotopic to $[h_n(x),h_n(y)]$. We
parametrise the arcs $[x,y]$, $[h_n(x),h_n(y)]$, $k$ and $k'$ by
their arc lengths.

By \cite[Lemma 4.2.10]{ceg}, for sufficiently large $n$, we have \linebreak
$d([x,y](t),k(t))\leq\eps'_n$
and $d([h_n(x),h_n(y)](t),k'(t))\leq\eps'_n$ for any $t$ with $\eps'_n\longrightarrow
0$.
 It follows that
$\displaystyle 1\leq\frac{l([x,y])}{l(k)}\leq\frac{l(k)+\eps'_n}{l(k)}\leq
1+\epsilon_n'$ and that $\displaystyle 1\leq\frac{l([h_n(x),h_n(y)])}{l(k')}\leq
1+\epsilon_n'$ for sufficiently large $n$, where $l(.)$ denotes the
length in $\bh^3/\rho_n(\pi_1(M))$. Therefore, we have
$\displaystyle d([x,y](\frac{l([x,y])}{l(k)}t),k(t))\leq 2\eps'_n$ for
sufficiently large $n$. For the same reason, we have also
$$d([h_n(x),h_n(y)](\frac{l([h_n(x),h_n(y)])}{l(k')} t),k'(t))\leq
2\eps'_n.$$

By the property (iii), the simple closed curve $\zeta_x\cup k\cup
\zeta_y\cup k'$ bounds a disc. Since $k$ and $k'$ are geodesic segments, the
function $d(k(t),k'(\frac{l(k')}{l(k)} t))$ is convex. Therefore we have
$d(k(t),k'(\frac{l(k')}{l(k)} t))\leq 5\eps$ for any $t$ since $d(x,h_n(x))\leq 5\eps$ and $d(y,h_n(y))\leq 5\eps$.

We define $g_n:[x,y]\rightarrow [h_n(x),h_n(y)]$ by the following
formula \linebreak
$g_n([x,y](\frac{l([x,y])}{l(k)}t))=[h_n(x),h_n(y)](\frac{l([h_n(x),h_n(y)])}{l(k)}t)$.
Setting $z=[x,y](\frac{l([x,y])}{l(k)}t)$, the distance
$d(z,g_n(z))$ is less than the following quantity
$$d([x,y](\frac{l([x,y])}{l(k)}t),k(t))+d(k(t),k'(\frac{l(k')}{l(k)}t))$$$$
+d(k'(\frac{l(k')}{l(k)}t,[h_n(x),h_n(y)](\frac{l([h_n(x),h_n(y)])}{l(k)}t)).$$ Then we get $d(z,g_n(z))\leq 5\eps+4\eps'_n$.
Now we conclude by
taking $N$ such that $4\eps'_n\leq \eps$ for $n\geq N$ and by
changing $A_n$ by a homotopy so that the geodesic segment
$\zeta_n(z)$ connecting $z$ to $g_n(z)$ lies in $A_n$ for any $z$ in
$\bar s$.
\end{proof}

Now we return to the proof of Lemma \ref{band}.
We remove from $\bar s$ its two extremal edges so that $g_n$ is defined on the
entire $\bar s$ (and we say that $\bar s$ originally had $j+2$ edges so that it now has $j$ edges).
By the equation (\ref{limit1}), there is a constant $R'$ such that for
any branch $\hat b$ of $\hat\tau^1$, we have $l(\tilde f_n(\hat
b))\leq R'\eps_n^{-1}$. By Lemma \ref{anne} and  since $\delta_n
\longrightarrow 0$ as $n \rightarrow \infty$, we have
$l_{\rho_n}(g_n(\bar s))\leq l_{\rho_n}(\bar s)+13\eps\leq j
R'\eps_n^{-1}+13\eps$. Therefore if we let $j'_n$ be the number of edges
contained in $g_n(\bar s)$, we have
\begin{eqnarray}\label{large}
j'_n\leq j (\frac{R'}{R})+1
\end{eqnarray}
 for large $n$, where $R$ is the constant in (b) of Lemma \ref{rtt}.

Since $\tau^1$ has only finitely many branches, once $j$ is fixed, there are only finitely
many possibilities for the oriented train routes of $\bar s$ and of $g_n(\bar s)$. Thus we
can find $\Theta(p_n)$-many disjoint segments $\{\bar  s_{1,n},...\}$ with the properties
(i), (ii) and (iii), having the same oriented train route, such that the segments $g_n(\bar s)$ also
have the same oriented train route for all $\bar s \in \{\bar s_{1,n}, \dots \}$.

Now we fix $j$  and extract a subsequence so that the oriented train routes of $\bar s_{i,n}$ and $g_n(\bar s_{i,n})$ do not depend on $n$.\\

By the same arguments, we can construct segments containing $j+1$ edges, $j+2$ edges and so on. To avoid  difficulties in the next part of the proof, we want to ensure that the segment containing $j$ edges which we have chosen is contained in the segment with $j+1$ edges as its first $j$ segments.
For that, to each $\bar s_{i,n}$, we add the edge of $\bar c_n$ which is adjacent to the last (with respect
to the orientation of $\bar s_{i,n}$) edge of $\bar s_{i,n}$.
 In this family of segments, take a
maximal family of disjoint segments lying in $\bar c_n^1$, which we denote by
$\{\bar s^+_{1,n},\bar s^+_{2,n},...\}$.
By Claim \ref{order of number} and the argument after that, there are only $o(p_n)$ segments among $\{\bar s_{1,n}, \dots \}$ for which the added edge lies outside of $\bar c_n^1$.
Therefore  the number of segments in $\{\bar s^+_{1,n}, \dots \}$ is $\Theta(p_n)$. It follows from the arguments we used for the segments $\bar s_{i,n}$ that among the segments $\bar s^+_{i,n}$ there are $\Theta(p_n)$ segments which have the properties (i), (ii) and (iii). The proof of Lemma \ref{anne} applies to these  $\Theta(p_n)$ segments,
yielding an homeomorphism $g_n$ which can be chosen to coincide with the one
defined on the segments $\bar s_{i,n}$ when restricted to them.

From this last family, we take $\Theta(p_n)$ segments $\bar s^+_{i,n}$ such that the oriented
train routes of $\bar s^+_{i,n}$ and $g_n(\bar s^+_{i,n})$ do not depend on $i$. Then we extract a
subsequence (with respect to $n$) such that for sufficiently large $n$, the oriented train
routes of $\bar s^+_{i,n}$ and of $g_n(\bar s^+_{i,n})$ do not depend on $n$.

By doing this argument recursively, increasing $j$ one by one, we complete the proof of
Lemma \ref{band}.
\end{proof}


\section{Commuting elements}    \label{annuli}

In this section, we shall use the results of the preceding sections to construct a sequence of discs which cross $\lambda$ less and less. This way, we shall obtain a homoclinic simple geodesic which does not cross $\lambda$.

Using Lemma \ref{band}, it is easy to observe that for $j$ (the
number of edges that $s_n$ contains) large enough, $s_n$ and $s'_n$
come back simultaneously to a branch of $\tau_1$. Still, the fact
that the edges of $\bar c_n^1$ are very long makes it difficult to
turn this observation into an actual construction. Instead we shall
use Lemma \ref{band} to construct two sequences $a_n, a'_n$ of
elements of $\pi_1(\partial M)$ which have nearly the same actions
on some rather large part of $\tilde f_n(\hat\tau_1)$. By results of
Kapovich (\cite{kapo}) this implies that the images of $a_n$ and
$a'_n$ in $\pi_1(M)$ commute and therefore correspond to an annulus.
Then we shall construct discs as we want by adding some waves.

Recall that given an oriented train
route $(bo(i))_{i\in\bn}$, we denote by $(b(i))_{i\in\bn}$ the corresponding
non-oriented train route. Since, by Lemma \ref{band}, for any $t$, the simple closed curve
$c_n$ goes through the oriented train routes $(bo(i))_{i\leq t}$ and
$(bo'(i))_{i\leq t}$ in $\tau^1$ for $n$ large enough,  there are two half-leaves
$l_+$ and $l'_+$ of the realised part of $L_{rec}$ whose oriented train routes are
$(bo(i))_{i\in\bn}$ and $(bo'(i))_{i\in \bn}$ respectively.

Let $e\in\bn$ be a natural number which we shall specify later, and let
$\varphi:\bn\rightarrow\bn$ be an increasing function such that:
\begin{enumerate}
\item[($\phi$1)] $bo(\varphi(i)+j)=bo(\varphi(0)+j)$ for any $i\in\bn$ and any $0\leq j\leq e$, i.e. the same train route
$(bo(\varphi(0)+j))_{0\leq j\leq e}$ is repeated starting from each $\varphi(i)$.
\item[($\phi$2)]
$bo'(V(\varphi(i))+j)=bo'(V(\varphi(0))+j)$ for any $i\in\bn$ and any $0\leq j\leq
 V(e)$, where $V$ is the function which appeared in Lemma \ref{band}.
\item[($\phi$3)] Suppose that $l_+$ (resp. $l'_+$) is not a closed curve, and let $k_{0,i}$ be a sub-arc
of $l_+$ with train route $(b(j))_{\varphi(0)\leq j<\varphi(i)}$ (resp. $k'_{0,i}$ the arc
of $l'_+$ with train route $(b'(i))_{V(\varphi(0))\leq i<V(\varphi(i))})$. The
two endpoints
of $k_{0,i}$ (resp $k'_{0,i}$) lie in the same switch of $\tau^1$ and the sequences
$(\partial
k_{0,i})$ (resp. $\partial k'_{0,i}$) converges to a single
point with respect to the Hausdorff topology as $i \longrightarrow \infty$.
\end{enumerate}

\begin{figure}
\psfrag{a}{$k_{0,1}$} \psfrag{b}{$bo(\varphi(1))$}
\psfrag{c}{$bo(\varphi(0))$} \psfrag{e}{$bo(\varphi(1)+e)$}
\psfrag{f}{$bo(\varphi(0)+e)$}
\includegraphics[width=.6\linewidth]{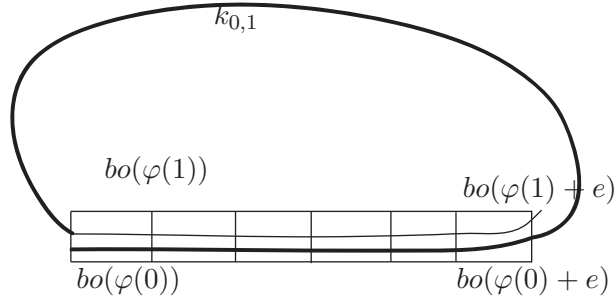}
\caption{Repeating the train route $(bo(\varphi(0)+j))_{0\leq j\leq e}$.}
\label{return}
\end{figure}

The existence of such a function $\varphi$ follows from the fact that the number of
branches and switches of $\tau^1$ is finite. See Figure \ref{return}.

 We define $\psi:\bn\to\bn$ as follows: $\psi(i)$ is the largest integer such that $\varphi(\psi(i))\leq i$. Since $\varphi$ is not surjective,
 it does not have an inverse function, and we use $\psi$ instead. Set $\psi_n=\psi(T(n)-e)$ for the function $T$ which appeared in Lemma \ref{band}.
 Note that if we forget the condition $(\phi 3)$, then $\psi_n+1$ is the number of times that $s_n$ comes back
 to the oriented train route $(bo(k))_{\varphi(0)\leq k \leq e+\varphi(0)}$ and $s'_n$ comes back to the oriented
 train route $(bo'(k))_{V(\varphi(0))\leq k \leq V(\varphi(0))+V(e)}$ at the same time. Since $\varphi$ is increasing and $T$ is unbounded, $(\psi_n)$ is unbounded.

Let us denote by $b'$ the  branch $bo'(V(\varphi(0)))$. By Lemma
\ref{band}, for any $i\leq \psi_n$, we have $g_n\circ f_n\circ
v(\varphi(i))\in f_n(b'(V(\varphi(i)))=f_n(b')$. It follows that,
for any $n$ large, there are $i_n,j_n\leq\psi_n$ such that the
distance between $g_n\circ f_n\circ v(\varphi(i_n))$ and $g_n\circ
f_n\circ v(\varphi(j_n))$ measured on $f_n(b')$ is at most
$\frac{1}{\psi_n} l_{\rho_n}(f_n)(b')=o(\eps^{-1}_n)$. Allowing
$g_n\circ f_n\circ v(\varphi(i_n))$ and $g_n\circ f_n\circ
v(\varphi(j_n))$ to be a bit further from each other, we can still keep
their distance to be $o(\eps^{-1}_n)$ while assuming that
$j_n-i_n\longrightarrow\infty$.

We pick such sequences $(i_n)$ and $(j_n)$, i.e.\ for each $n$, we
take two indices $i_n<j_n\leq\psi_n$ such that
$j_n-i_n\longrightarrow\infty$ and  the distance between
$g_n\circ f_n\circ v(\varphi(i_n))$ and $g_n\circ f_n\circ
v(\varphi(j_n))$ measured on $f_n(b')$ is $o(\eps^{-1}_n)$. We
denote by $I_n\subset s_n$ the  segment between $v(\varphi(i_n))$
and $v(\varphi(j_n))$ and by $J_n$ be the sub-segment of $s_n$
consisting of the $e$ vertices following $v(\varphi(j_n))$.

Let $\tilde s_n\subset\bh^3$ be a lift of $f_n(s_n)$, and let
$\tilde v(\varphi(i_n))$, $\tilde v(\varphi(j_n))$, $\tilde I_n$ and
$\tilde J_n$ be lifts of $f_n \circ v(\varphi(i_n))$, $f_n \circ
v(\varphi(j_n))$, $f_n(I_n)$ and $f_n(J_n)$ respectively, lying on
$\tilde s_n$.
We lift the map $g_n$ to a map $\tilde g_n$ from
$\tilde s_n$ to a lift $\tilde s'_n$ of $f_n(s'_n)$.
Let $\rho_n(a_n)\in\rho_n(\pi_1(M))$ be the covering translation which
takes $\tilde f_n(\tilde b(\varphi(i_n))$ to $\tilde f_n(\tilde b(\varphi(j_n))$.
Since $bo(\varphi(i_n)+j)=bo(\varphi(j_n)+j)$ for all $j\leq e$, the
isometry $\rho_n(a_n)$ acts as a translation on $\tilde
I_n\cup\tilde J_n$.

\begin{figure}[hbtp]
\psfrag{a}{\small $\tilde v(\varphi(i_n))$} \psfrag{b}{\small$\tilde
v(\varphi(j_n))$} \psfrag{I}{\small $\tilde I_n$} \psfrag{J}{\small
$\tilde J_n$} \psfrag{s}{\small $\tilde s_n$} \psfrag{g}{\small
$\tilde g_n$} \psfrag{r}{\small $\rho_n(a_n)$} \psfrag{c}{\small
$\tilde g_n(\tilde v(\varphi(i_n)))$} \psfrag{d}{\small $\tilde
g_n(\tilde v(\varphi(j_n)))$} \psfrag{t}{\small $\tilde s_n'$}
\psfrag{p}{\small $\rho_n(a_n')$} \psfrag{i}{\small $\tilde I_n'$}
\psfrag{j}{\small $\tilde J_n'$}
\centerline{\includegraphics[width=8cm]{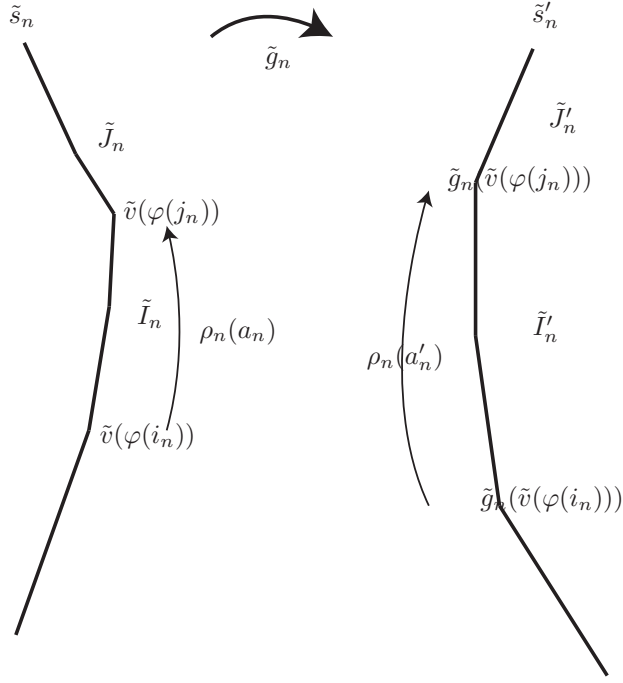}}
\caption{Commuting isometries} \label{translation}
\end{figure}

Let $\tilde I'_n\subset\tilde s'_n$ be the piecewise geodesic
segment between $\tilde g_n\circ\tilde v(\varphi(i_n))$ and
$\tilde g_n\circ\tilde v(\varphi(j_n))$, and let
$\tilde J'_n\subset\tilde s'_n$ be the segment between
$\tilde g_n\circ\tilde v(\varphi(i_n))$ and $\tilde
g_n\circ\tilde v(\varphi(j_n)+e)$. Let $\rho_n(a'_n)\in\rho_n(\pi_1(M))$ be the covering translation which
takes $\tilde f_n(\tilde b'(V(\varphi(i_n)))$ to $\tilde f_n(\tilde b'(V(\varphi(j_n)))$.
From the assumption that
$bo'(V(\varphi(i_n))+j)=bo'(V(\varphi(j_n))+j)$ for any $0 \leq j
\leq V(e)$, it follows that $\rho_n(a'_n)$ acts as a
translation on $\tilde I'_n\cup\tilde J'_n$.
See Figure \ref{translation}.

By our choice of $i_n,j_n$, we have  $$d(\tilde g_n\circ\tilde
v(\varphi(j_n)),\rho_n(a'_n)\circ\tilde g_n\circ\tilde v(\varphi(i_n))=o(\epsilon_n^{-1}).$$
 From this together with the facts that $\rho_n(a_n)$ acts as a translation on
$\tilde I_n\cup \tilde J_n$ and that $\rho_n(a'_n)$ acts as a translation on $\tilde
I'_n\cup\tilde J'_n$, we shall deduce the following claim:

\begin{cl}   \label{piti} For $R>0$, let $\cal{V}_R(\tilde J_n)$ be the
$R$-neighbourhood of $\tilde J_n$, then for any sequence of points $\tilde z_n\in
\cal{V}_R(\tilde
J_n)$, we have that $d(\tilde z_n,\rho_n(a^{-1}_n a'_n)(\tilde z_n))$ is $o(\epsilon_n^{-1})$.
\end{cl}
\begin{proof} It is sufficient to prove this claim for any sequence $(\tilde z_n)$ lying on $\tilde
J_n$.
Since $\rho_n(a_n^{-1})$ acts as a translation
on $\tilde I_n\cup\tilde J_n$, the point $\rho_n(a_n^{-1})(\tilde z_n)$ is the point
$\tilde z'_n\in\tilde I_n$ with
\begin{eqnarray}
d(\tilde z_n,\tilde v(\varphi(j_n)))=d(\tilde z'_n,\tilde
v(\varphi(i_n))),
\label{first equation}
\end{eqnarray}
 where $d$ denotes the  distance measured on  $\tilde I_n \cup \tilde J_n$.
%
The point $\tilde z''_n=\rho_n(a'_n)\circ\tilde g_n(\tilde z'_n) \in \tilde J'_n$ is  at the distance $d(\tilde
g_n(\tilde z'_n),\tilde g_n\circ\tilde v(\varphi(i_n))$ from $\rho_n(a'_n)\circ\tilde
g_n\circ\tilde v(\varphi(i_n))$ (measured on $\tilde I'_n\cup\tilde J'_n$).
As we saw above, $d(\tilde g_n\circ\tilde v(\varphi(j_n)),\rho_n(a'_n)\circ\tilde
g_n\circ\tilde v(\varphi(i_n)))=o(\epsilon_n^{-1}).$
Therefore we have $d(\tilde
z''_n, \tilde g_n\circ\tilde v(\varphi(j_n)))= d(\tilde g_n(\tilde z'_n),\tilde
g_n\circ\tilde v(\varphi(i_n)))+o(\epsilon_n^{-1})$.

By Lemma \ref{band}, $\tilde g_n$ moves each point within distance $6\epsilon$, and hence we get
the equality $$d(\tilde z''_n, \tilde g_n\circ\tilde v(\varphi(j_n)))=d(\tilde z'_n,\tilde
v(\varphi(i_n)))+o(\epsilon_n^{-1}).$$ Using equation \ref{first equation}, we get
$$d(\tilde z''_n, \tilde g_n\circ\tilde v(\varphi(j_n)))=d(\tilde z_n,\tilde
v(\varphi(j_n))+o(\epsilon_n^{-1}).$$ Since $\tilde J_n$ and $\tilde J_n'$ are within the Hausdorff distance $6\epsilon$ from each other, we have $$d(\tilde z''_n,\tilde
z_n)=o(\epsilon_n^{-1}).$$

By the triangle inequality
$$d(\rho_n( a'_n a^{-1}_n)(\tilde z_n),\tilde z_n)\leq d(\tilde z_n, \tilde z''_n)+ d(\tilde z''_n,\rho_n( a'_n a^{-1}_n)(\tilde z_n)).$$
and by $\tilde z''_n=\rho_n(a'_n)\circ\tilde
g_n\circ\rho_n(a_n^{-1})(\tilde z_n)$ we get
$$d(\tilde z''_n,\rho_n(a_n'a_n^{-1})\tilde z_n)=
 d(\rho_n((a'_n)^{-1})\tilde z''_n, \rho_n(a_n^{-1})\tilde z_n)$$$$=d(\tilde g_n\circ\rho_n(a_n^{-1})(\tilde z_n),\rho_n(a_n^{-1})
 (\tilde z_n))=o(\epsilon_n^{-1}).$$
Thus we finally get  $d(\rho_n(a'_n a^{-1}_n)(\tilde z_n),\tilde z_n)=o(\epsilon_n^{-1})$.
\end{proof}
We note that, as can be seen in the proof, the $o(\epsilon_n^{-1})$ is ``uniform'', namely there is a
sequence $\delta_n\longrightarrow 0$ independent of $(\tilde z_n)$ such that \linebreak$d(\rho_n(a'_n a^{-1}_n)(\tilde z_n),\tilde
z_n) \leq \delta_n \epsilon_n^{-1}$.

We shall use this claim to prove the following lemma.

\begin{lemma}  \label{ano} There is $N \in \bn$ such that for $n\geq N$, $\rho_n(a^{-1}_n
a'_n)=id$.
\end{lemma}

The main argument in the proof is  the following.

\begin{lemma}  \label{pitito} Let $[P_n,Q_n]\subset\bh^3$ be a sequence of geodesic
segments between $P_n$ and $Q_n$ such that $l([P_n,Q_n])$ is $\Theta(\eps_n^{-1})$ and
let $(\delta_n), (\delta'_n) \subset \pi_1(M)$ be two sequences such that the distances
$d(P_n,\rho_n(\delta_n)(P_n))$, $d(P_n,\rho_n(\delta'_n)(P_n))$,
$d(Q_n,\rho_n(\delta_n)(Q_n))$ and $d(Q_n,\rho_n(\delta'_n)(Q_n))$ are
all $o(\eps_n^{-1})$. Then there is $N$ such that for $n\geq N$,
$[\rho_n(\delta_n),\rho_n(\delta'_n)]=id$.
\end{lemma}
\begin{proof} This comes directly from the arguments which Kapovich used in
\cite[Theorem 10.24]{kapo} to prove that the minimal action of
$\pi_1(M)$ on $\br$-tree is small (see \cite[p.239]{kapo}).
\end{proof}

\begin{proof}[Proof of Lemma \ref{ano}]
Set $e=2p+1$ where $p$ is the number of the branches of $\tau^1$.
If we fix some
$n\in\bn$, by our choice of $e$, we can find two different integers $i_1$ and $i_2$ between $\varphi(j_n)$ and $\varphi(j_n)+e$ such
that $bo(i_1)=bo(i_2)$. Let $K_n\subset J_n$ be the segment of $s_n$ with train route
$(b(k))_{i_1\leq k\leq i_2}$, and let $\delta_n\in\pi_1(M)$ be an element which
takes
$\tilde v(i_1)$ to $\tilde v(i_2)$.
Since $i_1,i_2$ do not depend on $n$, the element $\delta_n$ does not depend on $n$ either. Let us denote it by $g$.
The isometry $\rho_n(g)$ acts as a translation
on the lift $\tilde K_n$ of $f_n(K_n)$ which lies in $\tilde s_n$.

Let $\tilde K_{n-}$ and $\tilde K_{n+}$ be the two extremal edges of $\tilde K_n$,
such that $\rho_n(g)$ takes $\tilde K_{n-}$ to $\tilde K_{n+}$. By Claim \ref{piti},
the translation length of  $\rho_n(a^{-1}_n a'_n)$ is $o(\epsilon_n^{-1})$ on $\tilde I_n\cup\tilde
J_n$. By the same arguments, the translation length of $\rho_n(a^{\prime-1}_n a_n)$ is also $o(\epsilon_n^{-1})$ on $\tilde I_n\cup\tilde J_n$. It follows that:
\begin{itemize}
\item for any sequence $\tilde z_n\in \cal{V}_R(\tilde K_{n+})$, $d(\tilde
z_n,[\rho_n(a^{-1}_n a'_n),\rho_n(g)](\tilde z_n))=o(\epsilon_n^{-1})$;\\
\item for any sequence $\tilde z_n\in \cal{V}_R(\tilde K_{n-})$, $d(\tilde
z_n,[\rho_n(g^{-1}),\rho_n(a^{-1}_n a'_n)](\tilde z_n))=o(\epsilon_n^{-1})$;\\
\item for any sequence $\tilde z_n\in \cal{V}_R(\tilde K_{n+})$, $d(\tilde
z_n,[\rho_n(g),\rho_n(a^{\prime-1}_n a_n)](\tilde z_n))=o(\epsilon_n^{-1})$.\\
\end{itemize}
 Since $\tilde K_{n+}$ and $\tilde K_{n-}$ are edges of $\bar c^1_n$,
their lengths are $\Theta(\eps_n^{-1})$.
Applying Lemma \ref{pitito} to the segments
$\tilde K_{n-}$ and $\tilde K_{n+}$, we see that for sufficiently large $n$, $[\rho_n(a^{-1}_n
a'_n),\rho_n(g)]$, $[\rho_n(g^{-1}),\rho_n(a^{-1}_n a'_n)]$ and
$[\rho_n(g),\rho_n(a^{\prime-1}_n a_n)]$ commute with $\rho_n(a^{-1}_n a'_n)$.
Therefore they belong to an elementary subgroup of $\rho_n(\pi_1(M))$. By \cite[p.239]{kapo}, it follows that the group generated by $\rho_n(a^{-1}_n a'_n)$ and
$\rho_n(g)$ is elementary.

Since $i_1\neq i_2$, the distance of translation of $\rho_n(g)$ is $\Theta(\eps^{-1}_n)$.
Since the group generated by
$\rho_n(a^{-1}_n a'_n)$ and $\rho_n(g)$ is elementary, there are $d\in\pi_1(M)$, $t,t_n\in\bn$
such that $g=d^t$ and $a^{-1}_n a'_n=d^{t_n}$. Since the  translation distance of $\rho_n(g)$ is $\Theta(\eps^{-1}_n)$,
so is the translation distance of $\rho_n(d)$. By Lemma \ref{piti} however, the translation distance of
$\rho_n(a^{-1}_n a'_n)$ is $o(\epsilon^{-1}_n)$. Therefore we have $\rho_n(a^{-1}_n
a'_n)=id$ for sufficiently large $n$.
\end{proof}

Now we are ready to construct the homoclinic geodesic which was announced at the beginning of the section.\\

\begin{lemma} \label{an}
Under the light assumptions, assume that $l_{\rho_n}(\lambda_n^*)$ is bounded, then there is a homoclinic simple
geodesic which does not cross $|\lambda|$.
\end{lemma}

\begin{proof}
Recall that, earlier in this section, we have picked a lift $\tilde s_n\subset\bh^3$ of $f_n(s_n)$, and lifts
$\tilde v(\varphi(i_n))$, $\tilde v(\varphi(j_n))$ of $f_n \circ v(\varphi(i_n))$ and $f_n \circ
v(\varphi(j_n))$ respectively, lying on
$\tilde s_n$. We have also lifted the map $g_n$ to a map $\tilde g_n$ from
$\tilde s_n$ to a lift $\tilde s'_n$ of $f_n(s'_n)$ and we have defined $a_n,a'_n\in \pi_1(M)$. 

To simplify the notations, set $\tilde x_n=\tilde v(\varphi(i_n))$ and $\tilde
y_n=\tilde v(\varphi(j_n))$. Recall that $\tilde
y_n=\rho_n(a_n)(\tilde x_n)$. By Lemma \ref{ano},
$\rho_n(a_n)(\tilde g_n(\tilde x_n))=\rho_n(a'_n)(\tilde g_n(\tilde
x_n))$, which implies that  $\rho_n(a_n)(\tilde g_n(\tilde x_n))$
lies on the same component of $\bar c_n^1$ as $\tilde g_n(\tilde
y_n)$. It follows that we can change $A_n$ so that
$\rho_n(a_n)(\tilde\zeta_n(\tilde x_n))$ (recall that
$\tilde\zeta_n(x_n)$ is a short segment in $\tilde A_n$ which joins
$\tilde x_n$ to $\tilde g_n(\tilde x_n)$) lies in $\tilde
A_n$. For the sake
of simplicity, we set $\tilde\zeta_n(\tilde
y_n)=\rho_n(a_n)(\tilde\zeta_n(\tilde x_n))$.

Denote by $\tilde k_n$ the segment of $\tilde c_n$ which is mapped by $\tilde f_n$ to an arc joining $\tilde x_n$ to $\tilde y_n$, namely $\partial \tilde f_n(\tilde k_n)=\{\tilde x_n,\tilde y_n\}$. Similarly denote by $\tilde k'_n$ the segment of $\tilde c_n$ which is mapped by $\tilde f_n$ to an arc joining $\tilde\zeta_n(\tilde x_n)$ to $\tilde\zeta_n(\tilde y_n)$. By construction, $\tilde k_n$ and $\tilde k'_n$ lie in lifts of $s_n$ and $s'_n$ respectively and goe through the train routes $(bo(i))_{\varphi(i_n)\leq
i<\varphi(j_n)}$ and $(bo'(i))_{V(\varphi(i_n))\leq
i<V(\varphi(j_n))}$ respectively.

The arcs $\tilde f_n(\tilde k_n)$, $\tilde f_n(\tilde k'_n)$, $\tilde\zeta_n(\tilde
x_n)$ and $\tilde\zeta_n(\tilde y_n)$ bound a disc $\tilde R_n$
in $\tilde A_n$ (by Lemma \ref{band}) and all of their endpoints lie
on the same component $\tilde f_n(\tilde c_n)$ of $\partial\tilde
A_n$. It follows that one component of $\tilde A_n-\tilde R_n$ is a
disc $\tilde D_n$. The boundary of this disc is the union of an arc $\tilde
d_n\subset \tilde f_n(\tilde c_n)$ and $\tilde\zeta_n(\tilde x_n)$
or $\tilde\zeta_n(\tilde y_n)$. See Figure \ref{nested}.

\begin{figure}[hbtp]
\psfrag{a}{$\tilde g_n(\tilde x_n)$}
\psfrag{b}{$\tilde g_n(\tilde y_n)$}
\psfrag{c}{$\tilde x_n$}
\psfrag{d}{$\tilde y_n$}
\psfrag{e}{$\tilde d_n$}
\psfrag{f}{$\bar c_n$}
\psfrag{g}{$\tilde c^*_n$}
\psfrag{D}{$\tilde D_n$}
\centerline{\includegraphics{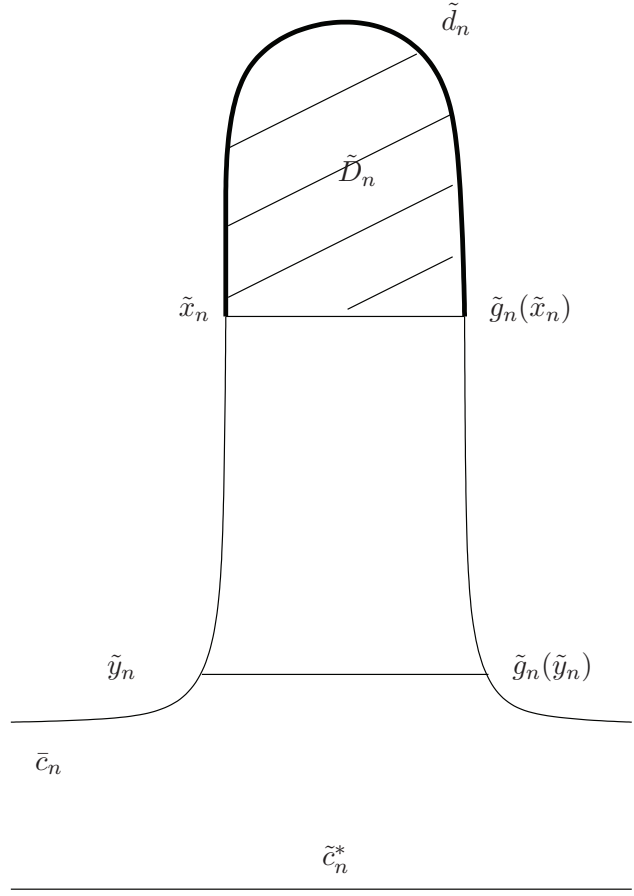}}
\caption{The innermost disc $\tilde D_n$}
\label{nested}
\end{figure}

\indent Let us fix a reference hyperbolic metric on $\partial M$. The endpoints of $k_n$ (resp. $k'_n$) are connected by an arc
$\kappa_n$ (resp. $\kappa'_n$) which lies on a switch of $\tau^1$. By Lemma \ref{ano}, the closed curves $k_n\cup\kappa_n$ and $k'_n\cup\kappa'_n$ are homotopic in $M$. In particular they bound an annulus $E_n$ in $M$.

Consider the projection $d_n\subset c_n$ of the arc which is mapped to $\tilde d_n$ by $\tilde f_n$. The closed curve $m_n=\kappa_n k_n d_n k_n^{\prime-1}\kappa^{\prime -1}_n d_n^{-1}$ (with appropriate choices of orientation) bounds a (possibly singular) disc which is the union of $E_n$ and two copies of the projection $D_n$ of $\tilde D_n$ (see Figure \ref{disc}).

\begin{figure}[hbtp]
\psfrag{a}{$\kappa_n$}
\psfrag{b}{$k_n$}
\psfrag{c}{$d_n$}
\psfrag{d}{$m_n$}
\centerline{\includegraphics{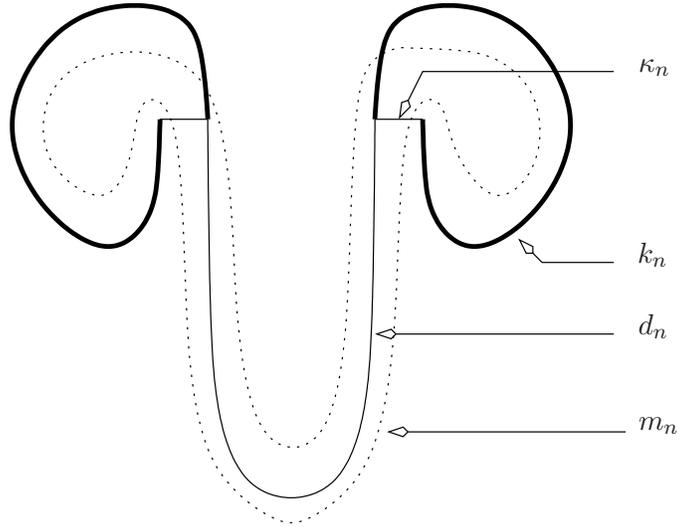}}
\caption{The curve $m_n$}
\label{disc}
\end{figure}

Since both $j_n-i_n$ and $i_n$ can be assumed to go to $\infty$, it follows from the property ($\phi$3)  that the length of $\kappa_n$
tends to $0$ (with respect to our reference metric on $\partial M$). It follows that $\int_{m_n}d\lambda_n\longrightarrow 0$. By the proof of the Loop Theorem (see \cite{hempel} for example), there is a meridian $m'_n$ such that $i(m'_n,\lambda_n)\longrightarrow 0$. By Casson's criterion (Theorem \ref{casson}) this concludes the proof of Lemma \ref{an}.
\end{proof}

\section{Conclusion}        \label{conclusion}
In this section, we shall complete the proofs of Theorem \ref{th&} and \ref{main}.

\begin{proof}[Proof of Theorem \ref{th&}]

Now we  complete the proof of Theorem \ref{th&}. Let $\rho_n$,
$\lambda_n$ and $\lambda$ be as in Theorem \ref{th&}. If no
subsequence of $(\rho_n)$ converges algebraically, it follows from
Lemma \ref{rtt} that $\rho_n$ and $\lambda_n$ satisfy the light
assumptions. Moreover $l_{\rho_n}(\lambda_n)$ is bounded by assumption. By Lemma \ref{an}, there is a
homoclinic geodesic which does not intersect
$|\lambda|$ transversely. By Lemma \ref{anho}, this contradicts the
assumption that $\lambda\in{\cal D}(M)$.
\end{proof}

We shall now deduce Theorem \ref{main} from Theorem \ref{th&}.
\setcounter{thm}{0}

\begin{thm}
Let $M$ be a compact irreducible atoroidal $3$-manifold with
boundary.
Let $(m_n)$ be a sequence in the Teichm{\"u}ller space $\mathcal{T}(
\partial M)$ which converges in the Thurston compactification to a projective
lamination $[\lambda]$ contained in $\cal{PD}(M)$. Let $q : \mathcal{T}(
\partial M)\rightarrow GF_0(M,P)$ be the Ahlfors-Bers map, and suppose that
$(\rho_n:\pi_1(M)\ra PSL(2,\bc))$ is a sequence of discrete
faithful representations corresponding to $(q(m_n))$. Then passing to a
subsequence, $(\rho_n)$ converges in $AH(N)$.
\end{thm}

\begin{proof} For a simple closed curve $c\subset
\partial M$, we denote by $l_{m_n}(c)$ the length of $c$ with respect to the
metric $m_n$ and by $l_{\rho_n}(c)$ the length of the closed geodesic of
$\bh^3/\rho_n(\pi_1(M))$ in the free homotopy class of $c$. By \cite[Theorem 2.2]{thuii} (see
also \cite{flp}), there is a sequence of simple closed curves
$c_n\subset\partial M$ whose projective classes converge to $[\lambda]$ in
$\cal{PML}(\partial
M)$ such that $\displaystyle\frac{l_{m_n}(c_n)}{l_{m_0}(c_n)}$ tends to $0$ as $n\rightarrow \infty$.

Using the following result of \cite{BC}, we shall get an upper bound for
the sequence $(l_{\rho_n}(\lambda_n))$.
\begin{thC}{\rm [Bridgeman-Canary]}
\label{B-C} For any $Q >0$, there is a constant $K>0$ depending only
on $Q$ with the following properties. Let $\Gamma$ be a finitely
generated Kleinian group without torsion such that the shortest
 length of the meridians with respect to the compatible hyperbolic structure on $\Omega_\Gamma/\Gamma$ is greater than $Q$. Let $C(\Gamma)$ be the convex
core of $\bh^3/\Gamma$, and consider the nearest point retraction
$r: \Omega_\Gamma/\Gamma \rightarrow
\partial C(\Gamma)$. Then $r$ is $K$-Lipschitz and has a homotopically inverse
$K$-Lipschitz map.
\end{thC}

Let us first verify that we are considering a situation where the
hypothesis of this theorem is fulfilled.
\begin{lemma}
\label{length grows}
There is a positive number $Q$ such that $l_{m_n}(d) \geq Q$ for any meridian $d$ of  $\partial M$.
\end{lemma}

\begin{proof} Assuming the contrary, we have a sequence of meridians $(d_n)$ such
that $(l_{m_n}(d_n))$ tends to $0$. Let us extract a subsequence so that
$(d_n)$ converges with respect to the Hausdorff topology to a geodesic lamination
$D\subset M$. By Casson's criterion, $D$ contains a homoclinic leaf. Since
$[\lambda]\in\cal{PD}(M)$, the lamination $D$ intersects the support of
$[\lambda]$ transversely (Lemma \ref{anho}). It follows that the sequence $i(\lambda,d_n)$ is bounded away
from $0$. This implies that the sequence $l_{m_n}(d_n)$ tends to $\infty$. Thus
we get a contradiction.
\end{proof}

It is clear that the length $l_{\rho_n}(c^*_n)$ is less than the length of any
curve in $
\partial C(\rho_n(\pi_1(M)))$ which is freely homotopic to $c_n$. Thus, applying
Theorem C, we see that $\displaystyle\frac{l_{\rho_n}(c^*_n)}{l_{m_0}(c_n)}$ tends to
$0$.

Let us denote by $\lambda_n$ the measured geodesic lamination obtained by
endowing $c_n$ with a Dirac measure whose weight is equal to
$l_{m_0}(c_n)^{-1}$. The sequence $\lambda_n$ converges in $\cal{ML}(\partial
M)$ to a measured geodesic lamination $\lambda$ which lies in the projective
class $[\lambda]$. Since $\displaystyle\frac{l_{\rho_n}(c^*_n)}{l_{m_0}(c_n)}$ tends to
$0$, we have $l_{\rho_n}(\lambda_n)\longrightarrow
0$. Since $\lambda$ lies in the projective class $[\lambda]\in{\cal PD}(M)$, the
measured geodesic lamination $\lambda$ lies in ${\cal D}(M)$. Applying Theorem
\ref{th&}, we see that a subsequence of $(\rho_n)$ converges algebraically.
\end{proof}

\section{Appendix}
For the sake of completeness, we shall give brief proofs of some propositions which we cited in previous sections.
\begin{prop1}{\rm [\cite{espoir}, \textsection $2$ paragraphs after Lemme 2.7]}  Let $\Sigma$ and $\Sigma'\subset
\partial_{\chi<0} M$ be two compact, connected, incompressible surfaces which are disjoint
or equal and do not contain any essential closed curve which can be homotoped into $
\partial_{\chi=0} M$. Let $\tilde\Sigma\subset
\partial\tilde M$ (resp. $\tilde\Sigma'$) be a connected component of the preimage of
$\Sigma$ (resp. $\Sigma'$) and let $\Gamma\subset\rho(\pi_1(M))$ (resp. $\Gamma'$) be the
stabiliser of $\tilde\Sigma$ (resp. $\Gamma'$).

Then $\overline{\tilde\Sigma} \cap\overline{\tilde\Sigma}'$ is either empty or
equal to the limit set of $\Gamma \cap\Gamma'$.

In the latter case, if $\Gamma \cap\Gamma'$ is not cyclic, then it is the fundamental
group of a (possibly twisted) $I$-bundle which is a connected component of a characteristic submanifold of
$(M,\Sigma \cup\Sigma')$. If $\Gamma \cap\Gamma'$ is cyclic, then it is a finite index
subgroup
of a solid torus which is a connected component of a characteristic submanifold of
$(M,\Sigma \cup\Sigma')$.
\end{prop1}
\begin{proof}Suppose that $\xi$ belongs to the limit set of $\Gamma\cap\Gamma'$.
Since we are considering only geometrically finite groups,
both $\Gamma$ and $\Gamma'$ are convex cocompact. Let $l$ be a geodesic ray from the origin $O\in \bh^3$ to $\xi$. Let $F$ be a fundamental
domain of the convex core of $\Gamma$ containing $O$, and $F'$ a fundamental
domain of the convex core of $\Gamma'$ containing $O$. Since $\Gamma$ and $\Gamma'$ are convex cocompact, the diameters of $F$ and $F'$
are bounded. Choose $g_n\in\Gamma, g_n'\in \Gamma'$ such that $g_n F\cap l\neq \emptyset, g_n'F'\cap l\neq\emptyset$ and
$g_nF\cap g_n'F' (\neq\emptyset) \longrightarrow \xi$. Since the diameters of $F$ and $F'$ are bounded, we have
$$d(g_nO, g_n'O)=d(O,g_n^{-1}g_n' O)\leq K$$ for all $n$ and a fixed $K$. By discreteness of $\rho(\pi_1(M))$, after passing to a subsequence,
$g_n^{-1}g_n'=g$ for all $n$.
Then $g_n^{-1}g_n'=g_m^{-1}g_m'$, i.e., $g_m g_n^{-1}=g_m'g_n'^{-1}\in \Gamma\cap\Gamma'$. Since $\Gamma$ and $\Gamma'$ are convex cocompact,
these elements are hyperbolic and the limit set of $\Gamma\cap\Gamma'$ contains at least two points.
Let $h\in\Gamma\cap\Gamma'$ be a hyperbolic element. Then the invariant geodesics  $\tilde c, \tilde c'$ of $h$ in $\tilde \Sigma,\tilde \Sigma'$
descend to closed geodesics $c$ and $c'$ in $\Sigma$ and $\Sigma'$. Hence $c$ and $c'$ bound an annulus $A$ (not necessarily embedded)
 which is not homotoped into $\partial M$. Hence $A$ must be contained in some characteristic submanifold of $(M,\Sigma\cup\Sigma')$.
If $\Sigma=\Sigma'$, then $c=c'$ and $A$ is a M\"{o}bius band.

Let $F=\overline{\tilde\Sigma} \cap\overline{\tilde\Sigma}'$ and $C,C'$ be the projections of the convex hulls of $F$ in $\tilde \Sigma$
and $\tilde \Sigma'$. Then $\pi_1(C)=\Gamma\cap\Gamma'$ and $C\cup C'$ is the boundary of an $I$-bundle $C\times I$ which is essential in
$(M,\Sigma \cup \Sigma')$. If $\Sigma=\Sigma'$, then it is a twisted $I$-bundle with the boundary $C$.
\end{proof}

Next we include some proofs from \cite{boches} for the reader's convenience.

\begin{prop}\label{ho}Suppose that  the action of $\pi_1(M)$ on an $\br$-tree $\cal T$ is minimal and small.
Let $S\subset \partial M$ be a compact compressible surface which has (possibly empty) geodesic boundary with respect to a hyperbolic metric on $\partial M$. Suppose that $\phi:\cal T_\mu \ra \cal T$ is a $\pi_1(M)$-equivariant morphism which folds only at
complementary regions.
If $\mu\subset S$ is in tight position with respect to a meridian $m\subset S$, then $|\mu|$ can be extended to a geodesic lamination which
contains a homoclinic leaf $h$ in $S$.
\end{prop}
\begin{proof}Fix a hyperbolic metric on $\partial M$. The universal cover $\tilde S\subset \bh^2$ is a convex subset of  $\bh^2$. Let $\bar m$ represent an element of $\pi_1(S)$
corresponding to the meridian $m$, which leaves invariant a lift $\tilde m\subset \tilde S$. By equivariance, for $x\in \tilde m$,
$$\phi(\pi(\bar m x))=\phi(\pi(x))$$ since $\bar m$ is trivial in $\pi_1(M)$, where $\pi:\tilde S\ra \cal T_\mu$ is an equivariant map. Since $\phi$ folds
$[\pi(x), \pi(\bar m x)]$ only finitely many times, one can find segments $\tilde I_1, \tilde I_2 \subset \tilde m$ such that
$\tilde I_1\cap\tilde I_2=y\in\tilde m$ and $\phi$ folds $\pi(\tilde I_1)$ and $\pi(\tilde I_2)$ along $\pi(y)$.
For $x\in m\cap \mu$, let $\mu_x^+$ denote a half leaf of $\mu$ starting from $x$ to a chosen positive direction. Then using Skora's
idea, Kleineidam and Souto showed \cite{boches} (Proposition 3) that there are $z_i\in I_i\cap \mu$ such that the lifts of $\mu_{z_1}^+$ and $\mu_{z_2}^+$ to
$\partial \tilde M$ have the same endpoints. Let $C$ be the complementary region of $\tilde \mu$ in $\tilde S$, which contains the folding
point $y$. Let $\tilde\mu_1,\tilde\mu_2$ be boundary leaves of $C$. Up to reversing the orientation, we can assume that $\tilde\mu_1^+,\tilde\mu_2^+$
are not asymptotic in $\tilde S$. Since the lifts of $\mu_{z_1}^+$ and $\mu_{z_2}^+$ to
$\partial \tilde M$ have the same end points, by shrinking the intervals, we can see that the projection of $\tilde\mu_1^+,\tilde\mu_2^+$
to $\partial \tilde M$ have the same endpoint. Let $l$ be the geodesic in $\tilde S$ joining the end points of $\tilde\mu_1^+,\tilde\mu_2^+$.
The projection of $l$ to $S$ becomes a homoclinic leaf disjoint from $\mu$.
\end{proof}

\begin{prop}\label{me}
Let $S$ be a compressible surface in $\partial M$, which contains a homoclinic leaf $h$. Then there is a sequence of meridians whose Hausdorff limit does not cross $h$.
\end{prop}
\begin{proof}
Since a homoclinic leaf cannot be contained in an incompressible surface by Lemma \ref{disque}, $S(\bar h)$ must contain a meridian $m$.

If $h$ contains  infinitely many homotopy classes of $m$-waves. Then there are $(x_i),(y_i)\subset \br$
such that $h(x_i),h(y_i)\in m$ and $h[x_i,y_i]$ are non-homotopic $m$-waves. Since $m$ is compact, after passing to a subsequence,
we may assume that $h(x_i)$ and $h(y_i)$ converge. Hence for any $\epsilon>0$, we can choose $i,j$ such that
the lengths of segments $[h(x_i),h(x_j)],[h(y_i),h(y_j)]\subset m$ are less than $\epsilon$.
Then $h[x_i,y_i]\cup h[x_j,y_j]\cup [h(x_i),h(x_j)]\cup[h(y_i),h(y_j)]$ is a meridian whose geodesic
representative lies nearby the homoclinic leaf.

If $h$ contains only finitely many homotopy classes of $m$-waves. then there is a meridian $m$ and two disjoint half-leaves $h^+$ and $h^-$ of $h$ such that $h^+$ and $h^-$ are in tight position with respect to $m$. Considering the intersections of $m$ and $h^+$ and $h^-$ respectively, one obtains a picture similar to Figure \ref{disc}, namely there is an arc $k\subset h^+$ and an arc $k'\subset h^-$ which nearly bounds an annulus and a wave between $k$ and $k'$. These arcs can be used to construct a sequence of meridians whose Hausdorff limit does not cross $h$ as explained in the proof of \cite[Proposition 1]{boches}.
\end{proof}

With more work, one could probably prove that $\bar h$ is a Hausdorff limit of meridians. On the other hand, in all the situations we have used Proposition \ref{me}, with only little changes, we could have replaced it with the following weaker result whose proof is easier.

\begin{lemma}\label{noncross}
Let $S$ be a compressible surface which contains a homoclinic leaf $h$, and let $\beta$ be a measured lamination which does not cross $h$. Then there is a sequence of meridians whose Hausdorff limit does not cross $\beta$.
\end{lemma}
\begin{proof}
Using cut-and-paste operation as in Claim \ref{limit}, we construct a sequence of meridians $m_i$ such that $i(m_i,\beta)\longrightarrow 0$. Start with a meridian $m\subset S(\bar h)$. Since $h$ is homoclinic, it contains an $m$-wave. Using this wave as in the proof of Claim \ref{limit}, we get a meridian $m_1$ such that $i(m_1,\beta)\leq\frac{1}{2} i(m,\beta)$. Then we do the same again on $m_1$. Repeating this, we get a sequence of meridians $m_i$ such that $i(m_i,\beta)\leq\frac{1}{2^i} i(m,\beta)$.
\end{proof}


\noindent     Inkang Kim\\
     School of Mathematics\\
    KIAS, Hoegiro 85, Dongdaemun-gu\\
     Seoul, 130-722, Korea\\
     e-mail: inkang\char`\@kias.re.kr\\

\noindent Cyril Lecuire\\
CNRS\\
Institut de Math\'{e}matiques de Toulouse\\
Universit\'{e} Paul Sabatier\\
118 route de Narbonne\\
31062 Toulouse Cedex 4\\
e-mail: lecuire\char`\@ math.ups-tlse.fr\\

\noindent
Ken'ichi Ohshika\\
Department of Mathematics\\
Graduate School of Science\\
Osaka University\\
Toyonaka, Osaka 560-0043, Japan\\
email: ohshika\char`\@ math.sci.osaka-u.ac.jp\\

\end{document}